# Two-point stress approximation: A simple and robust finite volume method for linearized (poro-)elasticity and Stokes flow

Jan Martin Nordbotten[1,2] & Eirik Keilegavlen[1]

1: Center for Modeling of Coupled Subsurface Dynamics,
Department of Mathematics, University of Bergen
Norway

2: Norwegian Research Center (NORCE)
Bergen, Norway

## Abstract

In this paper, we construct a simple and robust two-point finite volume discretization applicable to isotropic linearized elasticity, valid in also in the incompressible Stokes' limit. The discretization is based only on co-located, cell-centered variables, and has a minimal discretization stencil, using only the two neighboring cells to a face to calculate numerical stresses and fluxes. The discretization naturally couples to finite volume discretizations of flow, providing a stable discretization of poroelasticity.

We show well-posedness of a weak statement of the continuous formulation in appropriate Hilbert spaces, and identify the appropriate weighted norms for the problem. For the discrete approximations, we prove stability and convergence, both of which are robust in terms of the material parameters. Numerical experiments in 3D support the theoretical results, and provide additional insight into the practical performance of the discretization.

# 1. Introduction

There are strong similarities between the Laplace equation and the equilibrium equations of linearized elasticity, both superficially and in the sense of deep mathematical structures of relevance for the development numerical methods (see e.g. [1, 2]). Despite these similarities, while simple and robust discretizations of the Laplace equations are abundant and appear as canonical examples in any introductory textbook on discretizations of partial differential equations, the situation is much less straight-forward for linearized elasticity [3]. This paper addresses this gap by providing a simple and robust finite volume discretization applicable to elasticity and poroelasticity. As a corollary, it is also well suited as a co-located finite volume method for Stokes' equations.

The challenge of discretizing the equilibrium equations of linearized isotropic elasticity appears when considering the simplest form of these equations, namely [4]

$$\nabla \cdot (2\mu\epsilon(u) + \lambda(\nabla \cdot u)I) = f^u. \tag{1.1a}$$

Here $u$ is displacement, $\mu$ and $\lambda$ are Lamé parameters, and $\epsilon(u)$ is the linearized strain tensor typically expressed as:

$$\varepsilon(u) = \frac{1}{2}(\nabla u + \nabla u^T). \tag{1.1b}$$

When applying classical finite elements to discretize equation (1.1), the resulting method is not robust in the limit of incompressible materials, where $\lambda \to \infty$. Moreover, when coupling equation (1.1) to a porous or thermal material, instabilities appear when using lowest-order finite elements [5]. The situation is hardly better when considering finite volume or finite element methods, where the simplest stencils (e.g. five-point and seven-point stencils on logically Cartesian grids in 2D and 3D) cannot be constructed to provide a consistent discretization of $\epsilon(u)$.

As a result, a large literature has developed considering more complex discretization methods for linearized elasticity, of which we mention as examples stabilized finite element methods and mixed finite element methods (see an extensive survey in chapters 8 and 9 of [3], respectively), finite volume methods with larger stencils [6, 7, 8, 9], face-staggered variables [10] or nodal-staggered variables [11, 12].

In this paper we take a different approach. By considering an extended formulation of elasticity, we are able to provide a finite volume discretization that is both simple (in terms of the size of the discretization stencil) and robust (in terms of permissible grids and material parameters), as we make precise below.

Our development is based on combining two key insights. Firstly, we relax the symmetry of the stress tensor by extending the equations of elasticity to explicitly contain local rotations. This advantage of a "Hellinger-Reissner-type" formulation lies in the elegant connections to scalar-valued calculus [13, 14], which has proved to simplify the construction of numerical discretizations [14, 2, 15]. Secondly, by explicitly representing the solid mass density as an independent variable, we obtain robustness for incompressible materials (see [16] for a classical review of this idea). The solid mass density as an independent variable also allows for an algebraic coupling to e.g. poroelasticity [17]. We make the observation that the combination of these two perspectives results in a system of equations that can be manipulated in terms of a finite volume structure, with constitutive laws that only contain directional derivatives. This is the key observation that allows for the derivation of simple, consistent and robust finite volume methods.

To make the introductory comments more precise, we will consider an extended system of coupled partial differential equations modeling quasi-static poroelasticity on a domain $\Omega$. For introductory texts to mechanics and poroelasticity, we recommend the books of Temam and Miramville [4] and Coussy [18], respectively, which both follow a notation not dissimilar from ours.

As state variables, we consider a (macroscopic) deformation $u$ is seen from a Lagrangian perspective, a (microscale) solid rotation[1] $r_\circ$, a volumetric strain $\rho_\circ$ and the fluid pressure $w$. While the convention on variable names is modified to avoid duplicity, the notation below with respect to operators is standard in the analysis of elasticity. For a general review of the notation, and in particular the definition of the skew operators $S$ and $S^*$, see Appendix A1. The equations of linearized poroelasticity can then be summarized as (see e.g. [19, 18, 15]):

*Balance of linear momentum $\sigma$:*

$$\nabla \cdot \sigma = f^u. \tag{1.2a}$$

*Balance of angular momentum:*

$$-S\sigma = f_\circ^r. \tag{1.2b}$$

*Conservation of solid mass (equivalent to definition of volumetric strain):*

$$\nabla \cdot u - \rho_\circ = f^p. \tag{1.2c}$$

*Conservation of fluid mass flux $\chi$, in terms of a Biot modulus $\vartheta$ and compressibility $\eta_\circ$:*

$$\nabla \cdot \chi_\circ + \frac{\partial}{\partial t}(\vartheta \rho_\circ + \eta_\circ w) = f_\circ^w. \tag{1.2d}$$

For the sake of generality, we have retained right-hand side terms in all of equations (1.2), however in practice, one will often expect $f_\circ^r = f^p = 0$.

[1] The subscript ∘ denotes variables and functions that will be rescaled when writing the equations on a time-discrete conservation form in Section 1.1.

These balance equations are complemented by the linear constitutive laws, which in this work we restrict in generality to:

*Isotropic Hooke's law in terms of Lamé parameters $\mu$ and $\lambda$:*

$$\sigma = 2\mu(\nabla u + S^* r_\circ) + \lambda I \rho_\circ - \vartheta I w. \tag{1.3a}$$

*Darcy's law with permeability $\kappa_\circ$ and gravitational potential $h$:*

$$\chi_\circ = -\kappa_\circ \nabla(w + h). \tag{1.3b}$$

Note that to keep the presentation concise, we have simplified Darcy's law. In applications, the permeability will be factored into the permeability and fluid viscosity.

The governing equations above must be complemented by boundary conditions. We formulate the boundary conditions of Robin type, with one boundary condition for each constitutive law. We write this compactly by introducing on the boundary $\partial\Omega$ the non-negative definite invertible tensor-valued, possibly spatially varying, Robin weight $b^u$ and non-negative scalar-valued, possibly spatially varying, Robin weight $b^w$. Then for given elastic boundary data $g_u^u$ and $g_\sigma^u$, and corresponding fluid boundary data $g_w^w$ and $g_{\chi,\circ}^w$, we state boundary conditions in terms of the normal components of equations (1.3):

*Displacement-stress boundary conditions:*

$$b^u(\sigma \cdot n - g_\sigma^u) = 2\mu(g_u^u - u). \tag{1.4a}$$

*Fluid pressure-flux boundary conditions:*

$$b^w\big(\chi_\circ \cdot n - g_{\chi,\circ}^w\big) = \kappa_\circ\big(g_w^w - (w + h)\big). \tag{1.4b}$$

Here $n$ is the outer normal vector to the boundary $\partial\Omega$. Dirichlet boundary conditions are obtained for $b = 0$ and Neumann conditions are formally obtained in the limit $b^{-1} \to 0$.

An important observation is that the above general system of equations (1.2-1.4) is the common building block for several major field theories, which are often treated separately. We make this observation precise with the following remarks:

**Remark 1.1 [Linearized elasticity, $\vartheta = f_\circ^r = 0$]**: In absence of a right-hand side, equation (1.2b) simplifies to the condition that the stress tensor $\sigma$ is symmetric:

$$S\sigma = 0. \tag{1.5}$$

Eliminating further the porous structure by setting $\vartheta = 0$, equations (1.2a-1.2c) together with (1.3a) are the Hellinger-Reissner formulation of linearized elasticity [20], with the addition of the volumetric strain as an independent variable. This model can be simplified further, by noting that by inserting equation (1.3a) in equation (1.2b), we obtain $S(\nabla u + S^* r_\circ) = 0$. Since the symmetry operator satisfies $SS^* r_\circ = 2r_\circ$, Equation (1.2b) thus implies that

$$r_\circ = -\frac{1}{2} S\nabla u. \quad (1.6)$$

By eliminating the volumetric strain $\rho_\circ$ using equation (1.2c), and the rotation according to equation (1.6), the Hellinger-Reissner formulation of linearized elasticity simplifies to the form of elasticity expressed in terms of only displacement as a primary variable, given by equation (1.1). □

**Remark 1.2 [Stokes, $\vartheta = f_\circ^r = f^p = 0,\ \lambda \to \infty$]**: The steady-state Stokes equations are identical to the isotropic equations of linearized elasticity, with the interpretation of $u$ as a velocity [4]. For the Stokes equations, the incompressible limit is considered particularly important, i.e. $\lambda^{-1} \to 0$. When taking this limit, one changes variable from volumetric strain $\rho_\circ$ to pressure $p_\circ = \rho_\circ \lambda$. With this change of variables, the constitutive law (1.3a) takes the $\lambda$-independent form:

$$\sigma = 2\mu\varepsilon(u) + p_\circ I. \quad (1.7)$$

At the same time, the mass conservation equation becomes

$$\nabla \cdot u = \lambda^{-1} p_\circ. \quad (1.8)$$

The Stokes equations thus appear from equations (1.2a), (1.7) and (1.8) by formally considering the limit $\lambda^{-1} = 0$, and eliminating the stress variable:

$$\nabla \cdot \big(2\mu\varepsilon(u)\big) + \nabla p_\circ = f^u. \quad (1.9a)$$

$$\nabla \cdot u = 0. \quad (1.9b)$$

□

**Remark 1.3 [Poroelasticity]**: Equations (1.2-1.3), with $\vartheta > 0$, are the linearized equations for poroelasticity [18]. Indeed, the volumetric strain $\rho_\circ$ is the coupling variable from the mechanical deformation to fluid flow. On the other hand, the fluid pressure enters the constitutive law (1.3a) completely analogously to the solid pressure $p_\circ$. When considering a time-discrete formulation, the permeability $\kappa$ gets multiplied by the time-step, and from the discretization point of view it is therefore important to allow for the limit of $\kappa_\circ \to 0$ to avoid any time-step constraints [7]. We will elaborate on this model in Sections 1.1 and 4. □

**Remark 1.4 [Thermomechanics]**: By relabeling $w$ as enthalpy, we recognize that equation (1.2d) corresponds to conservation of energy while equation (1.3b) corresponds to Fourier's law [21]. The results of the present paper thus equally well apply to thermomechanics as poroelasticity. □

**Remark 1.5 [Porous media flow]**: While it is not the primary emphasis of this paper, it is clear that if $\vartheta = 0$, as discussed in Remark 1.1, the problem decomposes and one may still consider the equations of fluid mass conservation (1.2d) and Darcy's law (1.3b) as a model of flow in porous media [22]. Finite volume methods for this subsystem are well studied [23, 24], and the methods proposed in this paper are an extension of the so-called "two-point flux approximation", which is the standard method employed in the majority of academic and commercial codes for flow in geological porous media [25, 26, 27]. □

In view of the above, we make the following definition:

**Definition 1.6 [Robust]**: We reserve the word *robust* to imply any result that is valid for all limiting models discussed in Remarks 1.1 to 1.5, with degenerate parameters allowed as constrained by the Assumptions stated in Section 1.2. Moreover, a *robust discretization*, is robust independent of a (sufficiently small) grid parameter $\delta$. □

Our contributions in this paper can then be summarized as follows:

1. An algebraic reformulation of equations (1.2-1.4) into a conservation form well suited for finite volume discretization, given as equations (1.14-1.16) below.
2. *Robust* well-posedness theory for the extended equations of linearized micropolar elasticity given in equations (1.14-1.16).
3. A *robust* "two-point stress" finite volume approximation to equations (1.14-1.16), applicable to simplicial, Cartesian, and polyhedral grids. The proposed discretization has the following properties:
    a. The discrete system of equations is formulated only in terms of cell-centered primary variables.
    b. Stresses and fluxes across a face depend algebraically on only the two cells neighboring that face.
    c. The discrete system of equations arising from the finite volume structure has minimum matrix fill-in, in the sense that only cells sharing a face interact. As an example, this leads to 5-point and 7-point stencils in 2D and 3D on Cartesian grids.
4. Analysis proving the *robust* stability of the discretization on large class of grids.
5. Analysis proving the *robust* convergence of the discretization on grids possessing a "face orthogonal" property.
6. Numerical examples providing insight into the performance of the method in practice.

The rest of the paper is structured as follows. We continue this introduction with a derivation of the equations of elasticity on conservation form, and statement of the assumptions on material parameters used in this work. Section 2 provides the general finite volume discretization considered. Section 3 details the main contribution of this paper, which is the two-point stress approximation for the elastic sub-system. Section 4 provides the details for integrating the discretization of the elastic sub-system into a poromechanical discretization. Section 5 gives the theoretical justification for both the continuous and discrete form of the equations. Section 6 provides numerical verification of the proposed method and theoretical results. Brief concluding remarks are included in Section 7. The paper contains appendixes with detailed derivation of the method and an extended definition of the operators used in this work.

## 1.1 Model equations

While Equations (1.2-1.4) provide a context for the paper in line with the majority of existing literature, it is beneficial for the derivation of finite volume methods to consider an algebraically equivalent reformulation of these equations that directly allows for finite volume approximations.

In this section, and the remainder of the paper, we will consider the fluid mass conservation equation discretized in time by an implicit time-stepping method, and subsume the time-step into the permeability $\kappa = (dt)\kappa_\circ$ and flux $\chi = (dt)\chi_\circ$, while the previous time-step into the right-hand side $(dt)f_\circ^w \to f^w$. Moreover, we will for convenience omit the gravity term $h$ from Darcy's law.

To develop a simple finite volume discretization, it is necessary to rescale variables to *a rotation stress* $r$, a *displacement interpreted as a flux* $\upsilon$, an *asymmetry displacement tensor* $\tau$, and a *total pressure* $p$ according to:

$$r = 2\mu r_\circ, \qquad \upsilon = u, \qquad \tau = S^* u, \qquad \text{and} \qquad p = \lambda\rho_\circ - \vartheta w. \tag{1.10}$$

While the introduction of $\upsilon$ and $\tau$ may seem a bit arbitrary in the current context, their presence is justified by the developments below. Moreover, these formally equivalent notions of displacement appear naturally in the analysis of the continuous formulation (see e.g. [2]).

With these quantities introduced, we calculate from Equation (1.2d) that

$$\nabla \cdot \chi + \vartheta\rho_\circ + \eta_\circ w = \nabla \cdot \chi + \vartheta\lambda^{-1} p + \eta w. \tag{1.11}$$

Here we define the effective compressibility of the medium as

$$\eta = \eta_\circ + \lambda^{-1}\vartheta^2. \tag{1.12}$$

Furthermore, we note that Equation (1.2b) can be manipulated to conservation form, since:

$$-S\sigma = -S\big(2\mu(\nabla u + S^* r_\circ)\big) = -2\mu S\nabla u - SS^* 2\mu r_\circ = 2\mu\nabla \cdot \tau - 2r. \tag{1.13}$$

Here we have again used that $SS^* = 2$ and the fact that $S\nabla u = -\nabla \cdot S^* u = -\nabla \cdot \tau$, based on equations (1.10) and (A1.9).

In view of equations (1.11-1.13), equations (1.2) can be written on conservation form in terms of the new variables defined in equation (1.10):

$$\nabla \cdot \begin{pmatrix} \sigma \\ \tau \\ \upsilon \\ \chi \end{pmatrix} - \begin{pmatrix} 0 & & & \\ & \mu^{-1} & & \\ & & \lambda^{-1} & \vartheta\lambda^{-1} \\ & & -\vartheta\lambda^{-1} & -\eta \end{pmatrix} \begin{pmatrix} u \\ r \\ p \\ w \end{pmatrix} = \begin{pmatrix} f^u \\ f^r \\ f^p \\ f^w \end{pmatrix}. \tag{1.14}$$

In equation (1.14), we have rescaled the righthand side of equation (1.2b) $f^r = (2\mu)^{-1} f_\circ^r$. With the above change of variables, the constitutive laws given in equations (1.3) similarly take the following appealing form:

$$\begin{pmatrix} \sigma \\ \tau \\ \upsilon \\ \chi \end{pmatrix} = \begin{pmatrix} 2\mu\nabla & S^* & I & \\ S^* & & & \\ I & & & \\ & & & -\kappa\nabla \end{pmatrix} \begin{pmatrix} u \\ r \\ p \\ w \end{pmatrix}. \tag{1.15}$$

Finally, the boundary conditions (1.4) take the form:

$$\begin{pmatrix} b^u & \\ & b^w \end{pmatrix} \left( \begin{pmatrix} \sigma \\ \chi \end{pmatrix} \cdot n - \begin{pmatrix} g_\sigma^u \\ g_\chi^w \end{pmatrix} \right) = \begin{pmatrix} 2\mu(g_u^u - u) \\ \kappa(g_w^w - w) \end{pmatrix}. \tag{1.16}$$

As above, we note that $g_\chi^w = (dt) g_{\chi,\circ}^w$. We consider equations (1.14-1.16) as the poromechanical system on conservation form.

**Remark 1.7 [Properties of conservation form]**:
We emphasize the following key features of equations (1.14-1.16), that enable the construction of simple and robust finite volume discretizations:

1. Equations (1.14) are decomposed into a divergence term of generalized "stresses" and an algebraic mass matrix and coupling matrix.
2. The gradient terms in the constitutive laws in equation (1.15) only appear on the main diagonal, with solely scalar coefficients.
3. There is no interaction in the constitutive laws, equation (1.15), between the elastic "stress" variables $\sigma, \tau$ and $\upsilon$ and the fluid flux $\chi$. Thus, the derivation of appropriate finite volume schemes for the two subsystems completely decouple.

□

## 1.2 Assumptions on material parameters

Definition 1.6 sets the stage for the assumptions we will make in this manuscript. We make these assumptions precise as follows.

**Assumption 1.8 [Boundary conditions]**: We will not consider the case of pure Neumann boundary conditions, thus we assume that there are parts of the boundary with non-zero measure where $(b^u)^{-1} \neq 0$ and $(b^w)^{-1} \neq 0$. □

**Assumption 1.9 [Elastic subsystem, $\vartheta = 0$]**: When considering the elastic subsystem, we will make the assumptions that:

1. *Non-degenerate shear modulus:*

$$0 < \mu_- \leq \mu \leq \mu_+ < \infty. \tag{1.17}$$

2. *Lower bound on bulk modulus:*

$$0 < \lambda_- \leq \lambda. \tag{1.18}$$

□

**Assumption 1.10 [Porous subsystem, $\vartheta = 0$]**: When considering the porous media subsystem without elasticity, we make the assumption that:

1. *Non-degenerate permeability:*

$$0 < \kappa_- \leq \kappa \leq \kappa_+ < \infty. \tag{1.19}$$

□

**Assumption 1.11 [Poroelasticity, $\vartheta \geq 0$]**: When considering the coupled poromechanical system, we recall the definition of effective compressibility, $\eta = \eta_\circ + \vartheta^2\lambda^{-1}$, retain all assumptions listed in Assumption 1.9, and additionally assume:

1. *Non-negative and finite coupling to porous structure:*

$$0 \leq \vartheta\lambda^{-1} \lesssim 1. \tag{1.20}$$

2. *Non-negative permeability with upper bound:*

$$0 \leq \kappa \leq \kappa_+ < \infty. \tag{1.21}$$

3. *Non-negative fluid compressibility, with upper bound:*

$$0 \leq \eta_0 \leq \eta_+ < \infty. \tag{1.22}$$

4. *Non-degenerate fluid system:*

$$1 \lesssim \eta_\circ + \kappa. \tag{1.23}$$

□

The fourth point in Assumption 1.11 states that while the individual coefficients may all be degenerate according to points 1, 2, and 3 of the assumption, we cannot let all coefficients degenerate simultaneously. This would lead to a fluid system that was impermeable, incompressible, and also decoupled from the solid system. Such a system would clearly be unnatural to model as porous, since all terms on the left-hand side of equation (1.2d) evaluate to zero.

# 2. Finite volume discretization

We will apply the finite volume method to approximate equation (1.14-1.16). This requires us to first define the discrete representation of the computational domain.

## 2.1 Grid and discrete variables

We denote by $\mathcal{T}_\delta$ a family of finite non-overlapping partitions of $\Omega$ into *cells* $\omega_i \in \mathcal{T}_\delta$ associated with the index set $I_{\mathcal{T}_\delta}$.

**Definition 2.1 [Admissible grids]**: We consider the parameter $\delta > 0$ a representation of linear grid size, and require a grid to be *admissible* in the sense of:

1. Each cell is a *star shaped* polyhedron.
2. The cells have quasi-uniform size: $\delta \lesssim |\omega_i|^{\frac{1}{N}} \lesssim \delta$.
3. Each cell has a designated point $x_i \in \omega_i$ in the kernel of the cell, which we will refer to as the *cell center*.

We denote by $\mathcal{F}_\delta$ the $N-1$ dimensional *faces* of the partition, $\varsigma_k \in \mathcal{F}_\delta$ associated with the index set $I_{\mathcal{F}_\delta}$. We make the following requirement on faces:

4. Each face is a polygon.
5. Each face is the intersection of exactly two cells $\overline{\varsigma_k} = \overline{\omega_i} \cap \overline{\omega_j}$, or the intersection of a cell and the boundary $\overline{\varsigma_k} = \overline{\omega_i} \cap \overline{\partial\Omega}$.
6. The cell faces have quasi-uniform size: $\delta \lesssim |\varsigma_k|^{\frac{1}{N-1}} \lesssim \delta$.
7. Each face cell has a designated interior point $x_i \in \varsigma_k$, which we will refer to as the *face center*.

In view of point 5. above, and to facilitate the exposition of boundary conditions, we denote by $\mathcal{B}_\delta$ a nonoverlapping partition of the boundary $\partial\Omega$ into *boundary cells* $\varpi_i \in \mathcal{B}_\delta$, with index set $I_{\mathcal{B}_\delta}$, which will be treated as cells with degenerate width. We make the following requirements on boundary cells:

8. Each boundary cell is geometrically equivalent to a face, thus for any $\varpi_i$, there exists a $\varsigma_k$ such that $\varpi_i = \varsigma_k$, and the boundary cell center equals the face center $x_i = x_k$.
9. The index set of boundary cells is complementary to the index set of regular cells, thus $I_{\mathcal{T}_\delta} \cap I_{\mathcal{B}_\delta} = 0$.

□

Examples of admissible 2D grids are given in Figure 2.1, examples of admissible 3D grids can be seen in Figure 6.1.

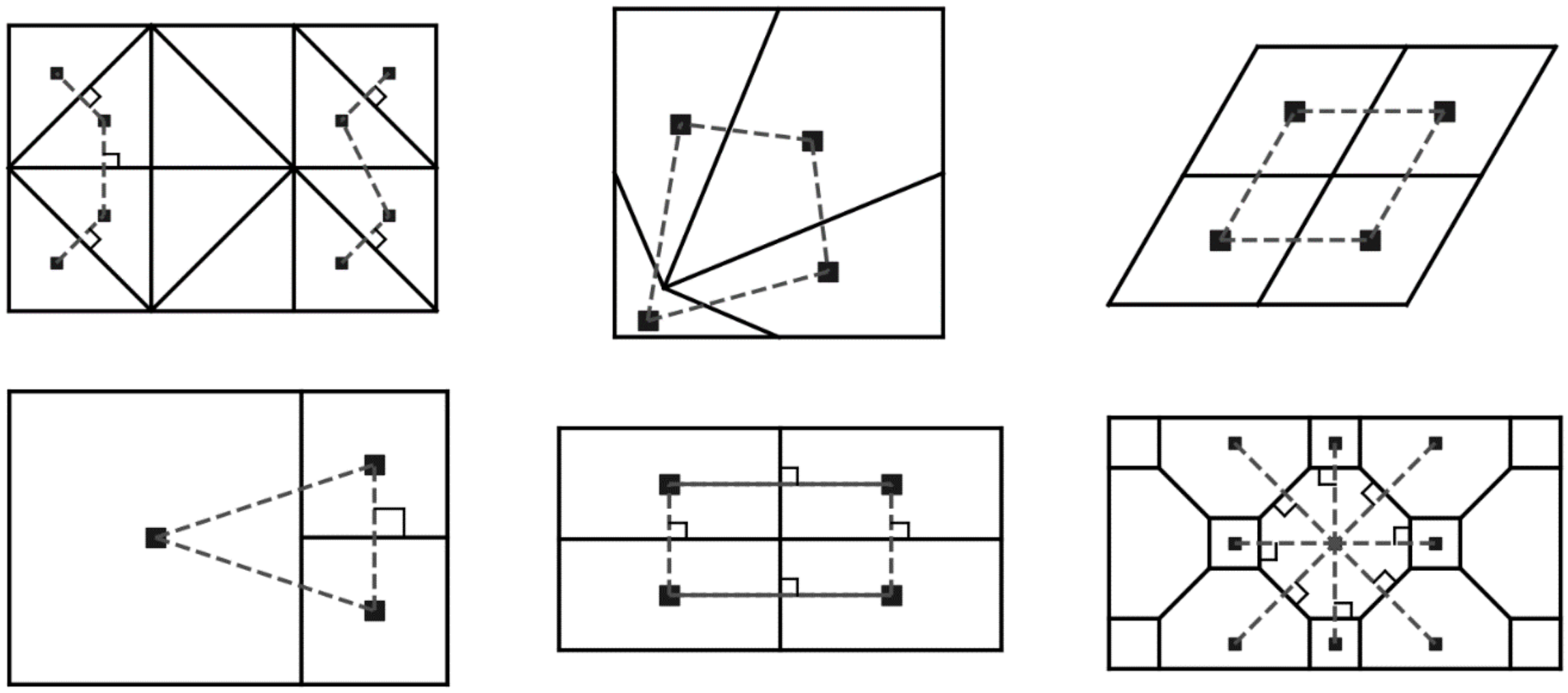

**Figure 2.1**: 2D illustrations of admissible meshes. We emphasize in particular A: Simplexes, B: Non-convex, C: Parallelogram, D: Local grid refinement: E: Rectangles and F: An Archimedean tiling. Face orthogonality is indicated by the right-angle symbol ∟

**Remark 2.2 [Curvilinear grids]**: Requirement 4 of Definition 2.1 states that the faces of grids be planar. This is for simplification of exposition, the extensions to non-planar grid faces (curved or piece-wise planar) are conceptually possible. □

The above construction of (internal and boundary) cells and faces is referred to as a *grid sequence*. We will often fix $\delta$, and consider a single given *grid*, and in this context we will frequently simply write $\mathcal{T}$, $\mathcal{F}$ and $\mathcal{B}$ to avoid carrying the subscript. We give the grid the following additional structure, which can be derived from the above.

For each cell $\omega_i$ and boundary cell $\varpi_i$ we denote the index set of neighbor faces $\mathcal{N}_i$. Note that the neighbor face of a boundary cell is given by point 8 above.

For each face $\varsigma_k$ we specify a normal vector $n_k$. We do not specify a convention on normal vectors for internal faces, but require that for boundary faces, $n_k$ is outward normal relative to the domain (and thus conforming with its continuous counterpart). The index set of neighbors is the adjoint of $\mathcal{N}$ and is denoted $\mathcal{N}_k^*$. Every face has two neighbors, and thus $\mathcal{N}_k^* = \{i, j\}$ and $|\mathcal{N}_k^*| = 2$ for all $k$. At each face, the vectors $n_{k,i}$ and $n_{k,j}$ point outwards from cell $i$ and $j$, respectively.

For all faces $\varsigma_k$, we denote the distances from the plane of the face to the centers of cells in $i \in \mathcal{N}_k^* \cap I_\mathcal{T}$ as

$$\delta_k^i \equiv (x_k - x_i) \cdot n_i. \tag{2.1}$$

Similarly, we denote the "distances" from the plane of the face to the centers of boundary cells $i \in \mathcal{N}_k^* \cap I_\mathcal{B}$ according to the boundary parameter in equations (1.4)

$$\delta_k^i \equiv b(x_i). \tag{2.2}$$

Note that when boundary conditions are uniform between the different variables, i.e. whenever $b^u = b^w I$ does not hold, equation (2.2) implies that the boundary distances $\delta_k^i$ will depend on what variable it is acting on in later expressions. We will not carry a separate notation for this, but highlight this when relevant. We further denote the (orthogonal projection) of the distances between the two adjacent cell centers $i, j \in \mathcal{N}_k^*$ as

$$\delta_k \equiv \sum_{i \in N_k^*} \delta_k^i. \tag{2.3}$$

The orientation implied by the normal vectors allows for the construction of an *incidence map*, $\Delta : (\mathcal{F}, \mathcal{B}) \to \mathcal{T}$. Concretely, we realize this map as the matrix $\Delta \in \mathbb{R}^{(|I_\mathcal{T}|+|I_\mathcal{B}|)\times|I_\mathcal{F}|}$ with entries $\Delta_{i,k} = n_{k,i} \cdot n_k$ if $i \in \mathcal{N}_k^*$ and $\Delta_{i,k} = 0$ if $i \notin \mathcal{N}_k^*$. The incidence map gives the connection between the definition of the grid structure and the finite volume structure, since it is clear that for any sufficiently smooth (vector) flux $\psi$, it holds that:

$$\int_{\partial\omega_i} \psi \cdot n \, dA = \sum_{k \in N_i} \int_{\varsigma_k} \psi \cdot n_{k,i} \, dA = \sum_{k \in I_{\mathcal{F}_\delta}} \Delta_{i,k} \int_{\varsigma_k} \psi \cdot n_k \, dA. \tag{2.4}$$

Complementary to the incidence map is a $\mu$-weighted averaging map $\Xi : \mathcal{T} \to \mathcal{F}$. Thus with $\mu_i$ denoting a cell-wise representation of the material parameter $\mu$, we realize $\Xi$ as the matrix $\mathbb{R}^{|I_\mathcal{F}|\times(|I_\mathcal{T}|+|I_\mathcal{B}|)}$ with entries

$$\Xi_{k,i} \equiv \frac{\mu_i}{\delta_k^i} \left( \sum_{j \in N_k^*} \frac{\mu_j}{\delta_k^j} \right)^{-1}. \tag{2.5}$$

Note that whenever $k$ is the index of a boundary face, then for $N_k^* = \{i, j\}$ we may without loss of generality consider $j \in I_\mathcal{B}$, and in the definition of $\Xi$, the boundary distance is always $b$. In view of the comment after equation (2.2), we note that unlike the incidence map, the *averaging map depends on the type of boundary condition, and therefore depends on what variable it acts on*. As an example, if a variable has a Dirichlet boundary condition, e.g. if $b^w = \delta_k^j = 0$, then we obtain from (2.5) that $\Xi_{k,i} = 0$ and $\Xi_{k,j} = 1$. Conversely, for a variable with a Neumann condition, e.g. if $(b^w)^{-1} = 0$ and thus formally $\delta_k^j = \infty$, then correspondingly $\Xi_{k,i} = 1$ and $\Xi_{k,j} = 0$. We will need both the averaging map and its complement, and define $\tilde{\Xi}_{k,i}$ as:

$$\tilde{\Xi}_{k,i} = \begin{cases} 1 - \Xi_{k,i} & \text{for} \quad k \in \mathcal{N}_i \\ 0 & \text{otherwise} \end{cases}. \tag{2.6}$$

For both the incidence map and the averaging map it will in the context of stating boundary conditions be useful to consider the restriction to internal and boundary cells. We denote these by subscripts, such that e.g. $\Delta_{\mathcal{F}} : \mathcal{F} \to \mathcal{T}$ and $\Delta_{\mathcal{B}} : \mathcal{B} \to \mathcal{T}$.

We close the grid description with one additional restrictions on grids, that will be necessary only in Section 5.3, when discussing the convergence of the discretization:

**Definition 2.3 [Face-orthogonal]:** We refer to a grid as *face-orthogonal* if $|n_{k,i} \times (x_k - x_i)| = 0$ for all faces $k$ and neighbor cells $i \in N_k^*$. □

We note that Definition 2.1 is very weak, and holds for essentially all grids composed of star-shaped polyhedra. On the other hand, face-orthogonality (which is required for a mesh to be "Admissible" in the sense of Eymard, Gallouët and Herbin [23]), is more restrictive, yet also holds for many classes of common grids provided the face and cell centers are appropriately defined. Examples include acute simplicial grids and many polyhedral grids, together with tensor product grids. Moreover, composite grids may often be face-orthogonal, such as any Archimedean tiling. These distinctions are illustrated in Figure 2.1.

## 2.2 Finite volume methods

In view of the definition of the divergence operator recalled in Appendix A1, as well as the definition of the grid and incidence matrix summarized by equation (2.4), we proceed to integrate equation (1.14) over each cell $\omega_i \in \mathcal{T}$. Considering therefore first the divergence term, we note that for any stress or flux $\psi \in \{\sigma, \tau, v, \chi\}$, it holds that:

$$\int_{\omega_i} \nabla \cdot \psi \, dV = \int_{\partial\omega_i} \psi \cdot n \, dA = \sum_{k \in I_{\mathcal{F}_\delta}} \Delta_{i,k} \int_{\varsigma_k} \psi \cdot n_k \, dA. \tag{2.7}$$

Based on equation (2.7), we are motivated to define the face traction (flux):

$$\psi_k \equiv \int_{\varsigma_k} \psi \cdot n_k \, dA. \tag{2.8}$$

$\psi_k$ is now a scalar or vector for each $k \in I_{\mathcal{F}}$, depending on whether $\psi$ is a vector or matrix. We combine all face values $\psi_k$ into the vector $\boldsymbol{\psi}$, such that

$$\sum_{k \in I_{\mathcal{F}_\delta}} \Delta_{i,k} \int_{\varsigma_k} \psi \cdot n_k \, dA = \sum_{k \in I_{\mathcal{F}_\delta}} \Delta_{i,k} \psi_k = \Delta \boldsymbol{\psi}. \tag{2.9}$$

Combining equation (2.7) and (2.9) now provides an exact expression of the divergence theorem, which we will refer to as the *finite volume structure*:

$$\int_{\omega_i} \nabla \cdot \psi \, dV = (\Delta \boldsymbol{\psi})_i. \tag{2.10}$$

Considering now the second term of equation (1.14), we note that for any primary variable $z \in \{u, r, p, w\}$ and material parameter $\gamma \in \{\mu, \lambda, \kappa, \eta, \vartheta\}$, it holds that

$$\int_{\omega_i} \gamma z \, dV = |\omega_i| (\pi_{\omega_i}\gamma)(\pi_{\omega_i} z) + \int_{\omega_i} (1-\pi_{\omega_i})\gamma \, (1-\pi_{\omega_i}) z \, dV, \tag{2.11}$$

where $\pi_{\omega_i}$ is the orthogonal projection onto the cell $\omega_i$. The deviatoric terms $(1-\pi_{\omega_i})$ are in general small (and will indeed be zero if the material parameter is constant on a cell). This motivates the definition of

$$\gamma_i = \pi_{\omega_i}\gamma, \qquad z_i = \pi_{\omega_i} z \quad \text{and} \quad f_i^z = \pi_{\omega_i} f^z. \tag{2.12}$$

For each of these cell-variables, parameters, and right-hand sides, we summarize them in vectors $\boldsymbol{\gamma}$, $\boldsymbol{z}$ and $\boldsymbol{f}^z$, in the same way we constructed $\boldsymbol{\psi}$ from $\psi_k$. Similarly, we denote by $|\boldsymbol{\omega}|$ the diagonal matrix with entries $|\omega_i|$ on the main diagonal.

With the above definitions, we obtain the finite volume structure for equation (1.14):

$$\Delta\boldsymbol{\psi} - |\boldsymbol{\omega}| \begin{pmatrix} 0 & & & \\ & \boldsymbol{\mu}^{-1} & & \\ & & \boldsymbol{\lambda}^{-1} & \boldsymbol{\vartheta}\boldsymbol{\lambda}^{-1} \\ & & \boldsymbol{\vartheta}\boldsymbol{\lambda}^{-1} & \boldsymbol{\eta} \end{pmatrix} \boldsymbol{z} \approx |\boldsymbol{\omega}| \boldsymbol{f}^z. \tag{2.13}$$

We emphasize that the expression is exact whenever the grid resolves the material coefficients, such that all material coefficients are constant within each cell $\omega_i \in \mathcal{T}$.

**Remark 2.4 [Finite volume methods]**: The structure presented in (2.13) defined *a finite volume method*. The particular finite volume method is thus distinguished by the construction of the numerical fluxes $\boldsymbol{\psi}$ from the cell-averaged quantities $\boldsymbol{z}$. □

# 3. Two-point approximation for linearized elasticity and Stokes

As follows from Remark 1.7, the elastic subsystem can be discretized independently of the flow subsystem due to the fact that there is no interaction between the two systems in the constitutive laws, equation (1.15).

In this section we therefore set $\vartheta = 0$, and consider the problem in the context of Assumption 1.9, deferring the discussion of the fluid pressure $w$ and the fluid flux $\chi$ until Section 4. For conciseness, we collect the remaining primary variables, "stresses" and the right-hand side in compound variables

$$y = (u, r, p)^T, \qquad \psi = (\sigma, \tau, v)^T, \qquad f = (f_u, f_r, f_p)^T. \tag{3.1}$$

The constitutive laws for the elastic subsystem can then be written as

$$\psi = \begin{pmatrix} 2\mu\nabla & S^* & I \\ S^* & & \\ I & & \end{pmatrix} y. \tag{3.2}$$

Similarly, the boundary condition for the elastic subsystem can be written as:

$$b^u(\sigma \cdot n - g_\sigma^u) = 2\mu(g_u^u - u). \quad (3.3)$$

## 3.1 Two-point stress approximations

From Section 2.2, we now know that the construction of an approximate representation of (3.2) is the key defining feature of a finite volume method. Locally, this construction will in general be a linear expression that can be expanded on the form:

$$\psi_k = \sum_{i \in N^*} T_k^i z_i + \sum_{i \in N^* NN_k^*} \hat{T}_k^i z_i + \cdots \quad (3.4)$$

Here, each $T_k^i$ and $\hat{T}_k^i$ is a linear map (matrix) from the primary variables $z_i$ to the space of stress variables $\psi_k$. The first right-hand side term contains contributions from the (at most) two neighbor cells of $\varsigma_k$, the second right-hand side term contains contributions from the (possibly many in 3D) neighbors-neighbor cells of $\varsigma_k$, and the dots indicate even more non-local connections. From this context we make the following definition:

**Definition 3.1 [Two-point stress approximation]:** The term "two-point stress approximation", abbreviated TPSA, refers explicitly to a numerical stress calculated based on an expression of the form (3.4), where the only non-zero coefficients are $T_k^i$. Thus, for any face $\varsigma_k \in \mathcal{T}$, the numerical stresses $\psi_k$ are approximated solely based on the two neighboring cells of $\varsigma_k$. □

To the authors knowledge, all existing consistent cell-centered finite volume methods for elasticity use "neighbors-neighbors" cells to approximate the traction $\sigma_k$ across a given face $k$, either explicitly (e.g. [6, 7, 8, 9]) or implicitly via staggered grids (e.g. [10] [11, 12]). This leads to rather large stencils, and more fill-in of the resulting system matrix than desirable. The use of neighbors-neighbors stress approximations is unavoidable if only the displacement is considered as a primary variable, since it is impossible to approximate rotations and volumetric stress based on only two points [28]. This observation implies that the extended set of variables introduced in Section 1.1 is strictly necessary for the construction of two-point stress approximations.

The choice of primary variables is therefore our key insight allowing for the current contribution. With this in mind, the TPSA stencil is a generalization of the construction of "two-point flux approximations" (TPFA) which are standard in commercial reservoir simulation. Building on theoretical analysis and experience with TPFA methods [23, 29], we expect *a priori* that the resulting discretization enjoys strong stability properties, but suffers from only being consistent on grids with high degree of symmetry. We will justify these expectations theoretically in Section 5, and numerically in Section 6.

The derivation of the TPSA coefficients is a bit technical and does not provide much insight. The derivation is therefore reported in full in Appendix A2; here we summarize the main results, which are sufficient for implementation and analysis.

The TPSA fluxes can be calculated in terms of the harmonic mean, $\bar{\mu}$, of the material coefficients $\mu$, together with a $\mu$-weighted distance $\delta^{\mu}$ (in these expressions, the shorthand $\delta_k^{-i} = \left(\delta_k^i\right)^{-1}$ is employed):

$$\bar{\mu}_k \equiv \delta_k \frac{\mu_i \delta_k^{-i} \mu_j \delta_k^{-j}}{\mu_i \delta_k^{-i} + \mu_j \delta_k^{-j}}, \quad \text{and} \quad \delta_k^{\mu} \equiv \frac{\left(\mu_i \delta_k^{-i} + \mu_j \delta_k^{-j}\right)^{-1}}{2}. \tag{3.5}$$

These combine with the difference operator and the $\mu$- and cell-face distance weighted mean values map defined in Section 2.1. Moreover, we introduce the short-hand notation for the rotation matrix $R_k^n$:

$$R_k^n = S^* n_k. \tag{3.6}$$

With this notation, Appendix A2 derives the numerical fluxes for an internal face $\varsigma_k$:

$$\begin{pmatrix} \sigma_k \\ \tau_k \\ \upsilon_k \end{pmatrix} = |\varsigma_k| \begin{pmatrix} -\delta_k^{-1} 2\bar{\mu}_k \Delta_k^* & -R_k^n \tilde{\Xi}_k & n_k \tilde{\Xi}_k \\ -R_k^n \Xi_k & & \\ n_k \Xi_k & & -\delta_k^{\mu} \Delta_k^* \end{pmatrix} \begin{pmatrix} \boldsymbol{u} \\ \boldsymbol{r} \\ \boldsymbol{p} \end{pmatrix}. \tag{3.7}$$

A comparison of Equations (3.2) and (3.7) reveals that they have the same overall structure, as should be expected. Of importance is the fact the material parameter $\mu$ appears as a harmonic means in the approximations for the normal derivative. The averaging terms $\Xi_k$ and $\tilde{\Xi}_k$ carry essentially the $\mu$-weighted values of the respective primary variables, and their combination is essential in preserving a quasi-symmetric structure of the discretization, as evidenced in Lemma 5.9 below. The new term arising in the lower-rightmost block of equation (3.7) can be interpreted as a discretization of $\delta_k^{-2} \mu^{-1} \nabla p$. This term arises naturally as part of the discretization, and has an important role in the stability of the method, as will be shown in Section 5.

Denoting finally as in Section 2.2 the composition of global operators from local operators by bold-face letters, we obtain the numerical fluxes for all faces as:

$$\begin{pmatrix} \boldsymbol{\sigma} \\ \boldsymbol{\tau} \\ \boldsymbol{\upsilon} \end{pmatrix} = |\boldsymbol{\varsigma}| \begin{pmatrix} -\boldsymbol{\delta}^{-1} 2\bar{\boldsymbol{\mu}} \Delta^* & -\boldsymbol{R}^{\boldsymbol{n}} \tilde{\Xi} & \boldsymbol{n} \tilde{\Xi} \\ -\boldsymbol{R}^{\boldsymbol{n}} \Xi & & \\ \boldsymbol{n} \Xi & & -\boldsymbol{\delta}^{\mu} \Delta^* \end{pmatrix} \begin{pmatrix} \boldsymbol{u} \\ \boldsymbol{r} \\ \boldsymbol{p} \end{pmatrix}. \tag{3.8}$$

In view of equation (A2.31) in Appendix A2.2, we note that equation (3.8) is also valid for homogeneous Dirichlet boundary conditions, due to the fact that the Robin weights $b^u$ and $b^r$ enter into the definition of the operators on the boundary faces as elaborated in Section 2.1. The appropriate inclusion of other boundary conditions is detailed in that appendix, and explicit expressions are provided in equations (A2.25-A2.28), but in the interest of space not included in the main text.

## 3.2 The finite volume TPSA discretization of linearized mechanics and Stokes

Equation (3.8) provides the numerical fluxes for the TPSA discretization of linearized elasticity. We combine these with the finite volume structure defined in Section 2.2 and equation (2.13), to obtain the global FV-TPSA discretization of linearized elasticity as:

$$\begin{pmatrix} -\Delta|\boldsymbol{\varsigma}|\boldsymbol{\delta}^{-1}2\overline{\boldsymbol{\mu}}\Delta^* & -\Delta|\boldsymbol{\varsigma}|R^{\boldsymbol{n}}\tilde{\Xi} & \Delta|\boldsymbol{\varsigma}|\boldsymbol{n}\tilde{\Xi} \\ -\Delta|\boldsymbol{\varsigma}|R^{\boldsymbol{n}}\Xi & -|\boldsymbol{\omega}|\boldsymbol{\mu}^{-1} & \\ \Delta|\boldsymbol{\varsigma}|\boldsymbol{n}\Xi & & -|\boldsymbol{\omega}|\boldsymbol{\lambda}^{-1} - \Delta|\boldsymbol{\varsigma}|\boldsymbol{\delta}^{\mu}\Delta^* \end{pmatrix}\begin{pmatrix} \boldsymbol{u} \\ \boldsymbol{r} \\ \boldsymbol{p} \end{pmatrix} = |\boldsymbol{\omega}|\begin{pmatrix} \boldsymbol{f}^u \\ \boldsymbol{f}^r \\ \boldsymbol{f}^p \end{pmatrix}. \tag{3.9}$$

We note that all expressions appearing in equation (3.9) are based on either differences or averages across a face, and that all terms combine at most a discrete divergence operator $\Delta$ with either an average $\Xi$ or a discrete directional derivative $\Delta^*$. Thus, the matrix in equation (3.9) has the minimum possible stencil for a finite volume discretization of a second-order elliptic PDE. As examples, for logically Cartesian grids in 2D and 3D, the proposed discretization results in 5-point and 7-point stencils, respectively.

We assert that the FV-TPSA system defined by equation (3.9) is robust in the sense of Definition 1.6 in the context of Assumptions 1.9. In particular, we assert setting $\lambda^{-1} = 0$ is a robust discretization of incompressible Stokes. We make this assertion precise, and prove it, in Section 5.2.1.

**Remark 3.2 [Multi-point stress]**: Since the term $|\boldsymbol{\omega}|\boldsymbol{\mu}^{-1}$ in equation (3.9) is diagonal, the rotation variable can be eliminated from the system according to

$$\boldsymbol{r} = -|\boldsymbol{\omega}|\boldsymbol{\mu}\Delta|\boldsymbol{\varsigma}|R^{\boldsymbol{n}}\Xi\boldsymbol{u}. \tag{3.10}$$

This comes at the cost of increasing the size of the discretization for the remaining variables $(\boldsymbol{u}, \boldsymbol{p})$. In particular, substituting (3.10) into (3.8) provides a multi-point stress approximation for the stress. Whether this elimination is advantageous will depend on the specifics of implementation and computer architecture. □

# 4. Two-point approximation for poroelasticity

Poroelasticity treats the problem of a porous solid, wherein the mechanical deformation is a result of both the solid stress and the fluid pressure. As pointed out in Remark 1.7, the model equations on conservation form given in Equations (1.14-1.16) decouple the discretization of poroelasticity into an elastic and Darcy subsystem.

The discretization of the elastic subsystem is presented in Section 3. For the purposes of this paper, any stable discretization of the flow system can be applied, as long as it is a finite volume method compatible with the grid structure of Section 2.1, and has cell-centered pressure variables. Examples of such methods are lowest-order mixed-finite elements [3], lowest-order control-volume finite element methods [30], multi-point flux approximation methods [31], and two-point flux approximation methods [25]. For all of the mentioned methods, we can construct a discrete flux approximation analogous to the expressing stated for stresses in Equation (3.4):

$$\chi_k = \sum_{i\in N^*} T_k^i w_i + \sum_{i\in N^* NN_k^*} \hat{T}_k^i w_i + \cdots \tag{4.1}$$

In the interest of the presentation being self-consistent and to have a concrete example for analysis, we include the two-point flux approximation (TPFA), which for homogeneous boundary conditions leads to the numerical flux [25, 23]:

$$\chi_k = |\varsigma_k| \delta_k^{-1} \bar{\kappa}_k \Delta_k^* \boldsymbol{w}. \tag{4.2}$$

Combining the numerical flux with the TPSA numerical stress given in equation (3.8) leads to:

$$\begin{pmatrix} \boldsymbol{\sigma} \\ \boldsymbol{\tau} \\ \boldsymbol{v} \\ \boldsymbol{\chi} \end{pmatrix} = |\boldsymbol{\varsigma}| \begin{pmatrix} -\boldsymbol{\delta}^{-1} 2\overline{\boldsymbol{\mu}}\Delta^* & -R^{\boldsymbol{n}}\tilde{\Xi} & \boldsymbol{n}\tilde{\Xi} & \\ -R^{\boldsymbol{n}}\Xi & & & \\ \boldsymbol{n}\Xi & & -\boldsymbol{\delta}^{\mu}\Delta^* & \\ & & & \boldsymbol{\delta}^{-1}\overline{\boldsymbol{\kappa}}\Delta^* \end{pmatrix} \begin{pmatrix} \boldsymbol{u} \\ \boldsymbol{r} \\ \boldsymbol{p} \\ \boldsymbol{w} \end{pmatrix}. \tag{4.3}$$

This numerical flux can be directly combined with the finite volume structure given in equation (2.13), leading to the poromechanical finite volume discretization:

$$\begin{pmatrix} -\Delta|\boldsymbol{\varsigma}|\boldsymbol{\delta}^{-1}2\overline{\boldsymbol{\mu}}\Delta^* & -\Delta|\boldsymbol{\varsigma}|R^{\boldsymbol{n}}\tilde{\Xi} & \Delta|\boldsymbol{\varsigma}|\boldsymbol{n}\tilde{\Xi} & \\ -\Delta|\boldsymbol{\varsigma}|R^{\boldsymbol{n}}\Xi & -|\boldsymbol{\omega}|\boldsymbol{\mu}^{-1} & & \\ \Delta|\boldsymbol{\varsigma}|\boldsymbol{n}\Xi & & -\Delta|\boldsymbol{\varsigma}|\boldsymbol{\delta}^{\mu}\Delta^* - |\boldsymbol{\omega}|\boldsymbol{\lambda}^{-1} & -|\boldsymbol{\omega}|\boldsymbol{\vartheta}\boldsymbol{\lambda}^{-1} \\ & & |\boldsymbol{\omega}|\boldsymbol{\vartheta}\boldsymbol{\lambda}^{-1} & \Delta|\boldsymbol{\varsigma}|\boldsymbol{\delta}^{-1}\overline{\boldsymbol{\kappa}}\Delta^* + |\boldsymbol{\omega}|\boldsymbol{\eta} \end{pmatrix} \begin{pmatrix} \boldsymbol{u} \\ \boldsymbol{r} \\ \boldsymbol{p} \\ \boldsymbol{w} \end{pmatrix} = |\boldsymbol{\omega}| \begin{pmatrix} \boldsymbol{f}^u \\ \boldsymbol{f}^r \\ \boldsymbol{f}^p \\ \boldsymbol{f}^w \end{pmatrix}. \tag{4.4}$$

# 5. Analysis

We structure the analysis in three parts. First, we will show the well-posedness of the continuous problem, in the sense of the four-field formulation used as the basis for the discretization. Secondly, we will provide the stability analysis for the FV-TPSA discretization, both for the mechanical subsystem and for the full poromechanical problem. Finally, we discuss consistency and convergence. In the interest of space, and since our main interest is in parameter and grid robustness, we will limit the exposition to homogeneous Dirichlet boundary conditions, the extension to more general boundary conditions being technical, but not conceptually challenging.

**Assumption 5.1 [Dirichlet boundary conditions]**: We make the blanket restriction that all calculations in Section 5 are made with homogeneous Dirichlet conditions for all variables, i.e. $b^u = b^w = g_u^u = g_\sigma^u = g_w^w = g_\chi^w = 0$. □

Throughout the section, we will collect as above the mechanical variables as $y = (u, r, p)$ and the poromechanical variables as $z = (y, w) = (u, r, p, w)$. We then make the notion of robust from Definition 1.6 and Assumptions 1.8-1.11 precise by introducing the solution norms:

$$\|y\|_*^2 = \|\nabla u\|_\mu^2 + \|u\|_\mu^2 + \|r\|_{\mu^{-1}}^2 + \|p\|_{\lambda^{-1}}^2 + |p|_{\mu^{-1}}^2 . \tag{5.1}$$

and

$$\|z\|_\circ^2 = \|y\|_*^2 + \|\nabla w\|_\kappa^2 + \|w\|^2 . \tag{5.2}$$

We remark that the norms in equations (5.1) and (5.2) represent explicitly the sufficient control and regularity we intuitively expect of the solution. Importantly, when $\kappa \to 0$, the regularity of the fluid pressure is reduced.

Throughout this section, we will use $\gamma$ (possibly with subscripts) to indicate the free parameter appearing from applying various inequalities.

## 5.1 Well-posedness of continuous formulation

We will show that the continuous formulation of poroelasticity as stated in Section 1.1 is well posed.

### 5.1.1 Elastic subsystem

We first consider the subsystem associated with linearized elasticity and Cosserat materials discussed in Section 3, and consider the norm stated in equation (5.1).

Define therefore the space of functions with bounded solution norm as:

$$Y = \{y \in L^2(\mathbb{R}^N) \times L^2(\mathbb{R}^N) \times L^2(\mathbb{R}) : \ \|y\|_*^2 < \infty\}. \tag{5.3}$$

We recognize that in the context of Assumption 1.9, this definition of $Y$ ensures that $u \in H^1$ and $r, p \in L^2$.

We now obtain the following weak formulation of the upper-left 3x3 block of the system of equations (1.14) and (1.15) by the usual approach of multiplying by a test function and integrating, the result being: Find $y \in Y$ such that:

$$-A_*(y, y') = (f^y, y') \qquad \text{for all} \qquad y' \in Y, \tag{5.4}$$

where:

$$\begin{aligned} A_*(y, y') = {} & (2\mu\nabla u, \nabla u') + (r, S\nabla u') + (p, \nabla \cdot u') + (S\nabla u, r') + \\ & +(\mu^{-1} r, r') - (\nabla \cdot u, p') + (\lambda^{-1} p, p') \end{aligned} \tag{5.5}$$

and

$$(f^y, y') = (f^u, u') + (f^r, r') + (f^p, p'). \tag{5.6}$$

Well-posedness of equation (5.4) relies on continuity and coercivity of the bilinear forms, which will be established in the two following Lemmas.

**Lemma 5.2 [Continuity of elastic subsystem, $\vartheta = 0$]**: Subject to Assumption 1.9, the bilinear form $A_*(y, y')$ is continuous with respect to the norm in equation (5.1), and $(f^y, y')$ is continuous in the $L^2$-norm.

*Proof*: The result is a direct application of the Cauchy-Schwarz inequality and the upper and lower bounds on $\mu$. □

We now recall the classical result for Stokes' equation.

**Lemma 5.3 [Stokes' inf-sup]**: For any function $p \in L^2(\mathbb{R})$ and for uniformly bounded weight $0 < \mu_- \leq \mu(x) \leq \mu_+ < \infty$, there exists a function $u_p \in H_0^1(\mathbb{R}^3)$ and $\beta > 0$ such that:

$$(\nabla \cdot u_p, p - \bar{p}) \geq \beta \left\| \nabla u_p \right\|_\mu |p|_{\mu^{-1}} \quad \text{and} \quad S\nabla u_p = 0. \tag{5.7}$$

Here $\bar{p}$ is the mean value of $p$, and $u_p$ is scaled such that

$$\left\| \nabla u_p \right\|_\mu = |p|_{\mu^{-1}}. \tag{5.8}$$

*Proof*: The result is shown in [32], without the spatial weight $\mu$. However, with the bounds stated in the Lemma, the weighted norms are equivalent to the unweighted norms, and the Lemma follows. □

This allows us to show coercivity, and we prove a slightly stronger result than needed in this section, as it will be useful when considering the coupled problem.

**Lemma 5.4 [Coercivity of elastic subsystem, $\vartheta = 0$]**: Subject to Assumption 1.9, the bilinear form $A_*(y, y')$ satisfies:

$$\inf_{y \in Y} \sup_{y' \in Y} \frac{A_*(y, y')}{\|y\|_* \|y'\|_*} \gtrsim 1. \tag{5.9}$$

*Proof*: We prove the Lemma by an explicit construction. Let any $y = (u, r, p) \in Y$ be given, and recall from Lemma 5.3 that there exists a $u_p$ satisfying equation (5.7). We then set $y' = (u + \alpha u_p, r, p)$, where $\alpha \in \mathbb{R}$ is a free parameter, and calculate:

$$A(y, y') = 2\|\nabla u\|_\mu^2 + 2(S\nabla u, r) + \|r\|_{\mu^{-1}}^2 + \|p\|_{\lambda^{-1}}^2$$

$$+\alpha\left(\left(2\mu\nabla u,\nabla u_p\right)+\left(r,S\nabla u_p\right)+\left(p,\nabla\cdot u_p\right)\right). \tag{5.10}$$

Before we proceed, we note the identities (recall that we only consider homogeneous boundary conditions in this analysis):

$$2\|\nabla u\|_\mu^2 = 2\|\varepsilon(u)\|_\mu^2 + \|S\nabla u\|_\mu^2 \quad \text{and} \quad \left(c,\nabla\cdot u_p\right)=0 \text{ for any } c\in\mathbb{R}. \tag{5.11}$$

Moreover, from Korn's inequality (A1.16), we know that

$$\|\varepsilon(u)\|_\mu^2 \geq C_K^2\|\nabla u\|_\mu^2\,. \tag{5.12}$$

Thus, by equation (5.11), (5.12) and Young's inequality

$$0\leq 2|(S\nabla u,r)| \leq \frac{1}{\gamma_1}\|S\nabla u\|_\mu^2+\gamma_1\|r\|_{\mu^{-1}}^2 = \frac{2}{\gamma_1}\left(\|\nabla u\|_\mu^2-\|\varepsilon(u)\|_\mu^2\right)+\gamma_1\|r\|_{\mu^{-1}}^2$$

$$\leq 2\frac{(1-C_K^2)}{\gamma_1}\|\nabla u\|_\mu^2+\gamma_1\|r\|_\mu^2. \tag{5.13}$$

Under the conditions of Assumption 5.1, it is straight-forward to calculate that $C_K^2 = 1/2$. Setting then $\gamma_1 = 1-\frac{C_K^2}{2}=\frac{3}{4}$, and using (5.11), (5.13) and (5.7), we obtain from (5.10) that:

$$A(y,y') \geq \frac{2}{3}\|\nabla u\|_\mu^2+\frac{1}{4}\|r\|_{\mu^{-1}}^2+\|p\|_{\lambda^{-1}}^2$$
$$+\alpha\left(-\|\nabla u\|_\mu|p|_{\mu^{-1}}+\beta|p|_{\mu^{-1}}^2\right). \tag{5.14}$$

Recognizing that we can use Young's inequality again with:

$$\|\nabla u\|_\mu|p|_{\mu^{-1}} \leq \frac{1}{\sqrt{\alpha}}\|\nabla u\|_\mu^2+\frac{\sqrt{\alpha}}{4}|p|_{\mu^{-1}}^2. \tag{5.15}$$

Then

$$A(y,y') \geq \left(\frac{2}{3}-\sqrt{\alpha}\right)\|\nabla u\|_\mu^2+\frac{1}{4}\|r\|_{\mu^{-1}}^2+\|p\|_{\lambda^{-1}}^2+\alpha\left(\beta-\frac{\sqrt{\alpha}}{4}\right)|p|_{\mu^{-1}}^2. \tag{5.16}$$

Choosing now any $\alpha$ sufficiently small, e.g.

$$\sqrt{\alpha}=\frac{1}{2}\min\left(\frac{2}{3},4\beta\right), \quad (5.17($$

the coefficients in front of all norms are positive, dependent only on constant in the Stokes' inf-sup, thus in light of the Poincaré inequality we obtain

$$A(y,y') \gtrsim \|y\|_*^2. \tag{5.18}$$

The Lemma follows since $\|y'\|_* \lesssim \|y\|_*$.

□

**Corollary 5.5 [Well-posedness of the elastic subsystem]**: Subject to Assumption 1.9, equation (5.4) is well-posed, with a weak solution $y \in Y$ satisfying

$$\|y\|_* \lesssim \|f^y\|. \tag{5.19}$$

*Proof*: The corollary follows Lemma 5.2, 5.4 and standard saddle-point theory, see Section 4.2.3 of [3] □.

### 5.1.2 Full poromechanical system

We consider now the full poromechanical system as stated in Section 1.1, in terms of the norm stated in equation (5.2).

Define first for the Darcy subsystem the space of functions with bounded $\kappa$-weighted norm as:

$$W = \{w \in L^2(\mathbb{R}) : \ \|\nabla w\|_\kappa^2 + \|w\|^2 < \infty\}. \tag{5.20}$$

We recognize that this definition of $W$ captures the intuitive expectation of how the regularity of $w$ depends on $\kappa$. We then define for the full poromechanical system the space of functions with bounded solution norm as

$$Z = Y \times W. \tag{5.21}$$

Note that the dependency on $\eta$ is not included for the fully coupled problem, which is a key result of the analysis.

We now obtain the following weak formulation of the system of equations (1.14-1.16) by the usual approach of multiplying by a test function and integrating, the result being: Find $z \in Z$ such that:

$$-A_\circ(z, z') = (f^z, z') \qquad \text{for all} \qquad z' \in Z, \tag{5.22}$$

where $z = (y, w) = (u, r, p, w)$ and:

$$A_\circ(z, z') = A_*(y, y') - (\vartheta\lambda^{-1}p, w') + (\vartheta\lambda^{-1}w, p') - (\kappa\nabla w, \nabla w') - (\eta^{-1}w, w'), \tag{5.23}$$

and

$$(f^z, z') = (f^u, u') + (f^r, r') + (f^p, p') + (f^w, w'). \tag{5.24}$$

Well-posedness of equation (5.22) again relies on continuity and coercivity of the bilinear forms, which will be established below.

**Lemma 5.6 [Continuity of poromechanical system]**: Subject to Assumptions 1.11, the bilinear forms $A_\circ(z,z')$ and $(f^z,z')$ are continuous with respect to the norm in equation (5.2).

*Proof*: The result is still a direct application of the Cauchy-Schwarz inequality and the stated bounds on material coefficients. □

We now introduce the following useful Lemma, that ensures that also the full poromechanical system is coercive. We state this in slightly more generality than what is needed at present, as the Lemma will be reused in Section 5.2.

**Lemma 5.7 [Coercivity of coupled poroelasticity]**: For $z=(y,w)\in Z$ and $z'\in Z$, consider a bilinear form:

$$A_\circ(z,z') = A_*(y,y') - B(y,w') + B(w,y') - C(w,w'). \qquad (5.25)$$

Then if Assumption 1.11 holds and:

1) Whenever Assumption 1.9 holds, the elastic sub-system $A_*$ has a coercivity estimate on the form
$$\inf_{y\in Y}\sup_{y'\in Y}\frac{A_*(y,y')-\|p\|_{\lambda^{-1}}\|p'\|_{\lambda^{-1}}}{\left(\|y\|_*-\|p\|_{\lambda^{-1}}^2\right)\left(\|y'\|_*-\|p'\|_{\lambda^{-1}}^2\right)}\gtrsim 1, \qquad (5.26)$$
2) Whenever Assumption 1.10 holds, the Darcy system $C$ has a coercivity estimate of the form:
$$\inf_{w\in W}\frac{C(w,w)-\|w\|_\eta^2}{\|\nabla w\|_\kappa^2}\gtrsim 1, \qquad (5.27)$$
3) The coupling terms are on the form (for $y=(u,r,p)$):
$$B(y,w') = (\vartheta\lambda^{-1}p, w'), \qquad (5.28)$$

Then the bilinear form $A_\circ(z,z')$ satisfies:

$$\inf_{z\in Z}\sup_{z'\in Z}\frac{A_\circ(z,z')}{\|z\|_\circ\|z'\|_\circ}\gtrsim 1. \qquad (5.29)$$

*Proof*: Let $z=(y,w)$ be given, and with respect to this $y$, let $y'=\left(u+\alpha u_p, r, p\right)$ be the vector defined in the proof of Lemma 5.4. We now set $z'=(y',-w)$ and calculate:

$$\begin{aligned} A_\circ(z,z') &= A_*(y,y') + B(y,w) + B(w,y') + C(w,w) \\ &\geq C_1\left(\|y\|_*^2 - \|p\|_{\lambda^{-1}}^2\right) + \|p\|_{\lambda^{-1}}^2 + 2(\vartheta\lambda^{-1}p,w) + C_2\|\nabla w\|_\kappa^2 + \|w\|_\eta^2, \end{aligned} \qquad (5.30)$$

where $C_1$ and $C_2$ are the hidden constants in inequalities (5.9) and (5.27). Then we once more use Young's inequality with $\gamma<1$ to obtain:

$$2(\vartheta\lambda^{-1}p, w) \leq \gamma\|p\|_{\lambda^{-1}}^2 + \frac{1}{\gamma}\|w\|_{\vartheta^2\lambda^{-1}}^2. \tag{5.31}$$

Combining equations (5.30), (5.31) and the Poincaré inequality (A1.15), we obtain

$$A_\circ(z, z') \geq C_1\left(\|y\|_*^2 - \|p\|_{\lambda^{-1}}^2\right) + (1-\gamma)\|p\|_{\lambda^{-1}}^2 + \frac{C_2}{2}\|\nabla w\|_\kappa^2 + \|w\|_{\eta-\gamma^{-1}\vartheta^2\lambda^{-1}+\frac{C_2 C_P}{2}\kappa}^2. \tag{5.32}$$

where $C_P$ is the Poincaré constant. In view of point 4 of Assumption 1.11, the weight on the $L^2$ term is bounded from below, in particular we can choose

$$\eta - \gamma^{-1}\vartheta^2\lambda^{-1} + \frac{C_2 C_P}{2}\kappa = \eta_0 + (1-\gamma^{-1})\vartheta^2\lambda^{-1} + \frac{C_2 C_P}{2}\kappa \gtrsim 1 + (1-\gamma^{-1})\vartheta^2\lambda^{-1}$$

Since by Assumption 1.11, we have that $\vartheta\lambda^{-1} \lesssim 1$, and we can choose $\gamma$ sufficiently close to 1 such that both

$$\eta - \gamma^{-1}\vartheta^2\lambda^{-1} + \frac{C_2 C_P}{2}\kappa \gtrsim 1 \qquad \text{and} \qquad (1-\gamma) > 0. \tag{5.33}$$

The lemma now follows from equation (5.32) and the definition of $\|z\|_\circ^2$. □

We summarize the above in the following Theorem.

**Theorem 5.8 [Well-posedness of coupled poroelasticity]**: Subject to Assumption 1.11, equation (5.22) is well-posed, with a weak solution $z \in Z$ satisfying

$$\|z\|_\circ \lesssim \|f^z\|.$$

*Proof*: The corollary follows Lemma 5.6, and 5.7 and standard saddle-point theory [3]. To show that the conditions of Lemma 5.7 hold, we remark that condition 1 follows from equation (5.18) in the proof of Lemma 5.4, and condition 3 holds by definition.

It thus only remains to show condition 2. This is however an immediate consequence of the definition of the bilinear form and the (weighted) norm. □

## 5.2. Well-posedness of FV-TPSA

The analysis of the well-posedness of the FV-TPSA discretization will follow the same general structure as for the continuous problem. This is in particular thanks to Lemma 5.9, which shows that the discretization indeed inherits many of the structural properties of the continuous problem.

In this section we will only discuss whether the FV-TPSA is well-posed, in the sense that for admissible grids (Definition 2.1) the solution is well-defined depending only on robust constants independent of the grid parameter $\delta$. The approximation properties of the method are discussed in Section 5.3.

For the discrete variables, we use spaces and norms analogous to those defined in equation (5.1-5.2). In particular, the vectors of discrete variables are denoted $\boldsymbol{y} = (\boldsymbol{u}, \boldsymbol{r}, \boldsymbol{p}) \in \boldsymbol{Y}$ and $\boldsymbol{z} = (\boldsymbol{y}, \boldsymbol{w}) = (\boldsymbol{u}, \boldsymbol{r}, \boldsymbol{p}, \boldsymbol{w}) \in \boldsymbol{Z}$, where

$$\boldsymbol{Y} = \mathbb{R}^{(N\times|I_\mathcal{T}|)} \times \mathbb{R}^{(N\times|I_\mathcal{T}|)} \times \mathbb{R}^{|I_\mathcal{T}|}, \qquad \boldsymbol{W} = \mathbb{R}^{|I_\mathcal{T}|} \tag{5.34}$$

and

$$\boldsymbol{Z} = \boldsymbol{Y} \times \boldsymbol{W} \tag{5.35}$$

These vectors are normed by:

$$\|\boldsymbol{y}\|_*^2 = \|\boldsymbol{\delta}^{-1}\Delta\boldsymbol{u}\|_{\bar{\mu}}^2 + \|\boldsymbol{u}\|_{\mu}^2 + \|\boldsymbol{r}\|_{\mu^{-1}}^2 + \|\boldsymbol{p}\|_{\lambda^{-1}}^2 + |\boldsymbol{p}|_{\mu^{-1}}^2 \tag{5.36}$$

and

$$\|\boldsymbol{z}\|_\circ^2 = \|\boldsymbol{y}\|_*^2 + \|\boldsymbol{\delta}^{-1}\Delta\boldsymbol{w}\|_{\bar{\kappa}}^2 + \|\boldsymbol{w}\|^2 \tag{5.37}$$

It should be clear from context that all norms involving expressions with $\Delta$ are face-norms, while the remaining norms are cell-norms (see equation (A1.12) for the distinction). Moreover, we emphasize that for all the face-norms, the weights appearing in the norms are always the harmonic average of the adjacent cell values.

As in the continuous case, the material weights appearing in the norm definitions (5.36) and (5.37) ensure the proper notion of regularity when variables degenerate. Notably, when $\kappa \to 0$, the regularity of $\boldsymbol{w}$ is reduced.

### 5.2.1 Elastic subsystem

In this section, we consider the FV-TPSA discretization stated in equation (3.9). By designating the discretization matrix of that equation by $\boldsymbol{A}_*$, the FV-TPSA discretization can be written compactly as:

$$\boldsymbol{A}_*\boldsymbol{y} = |\boldsymbol{\omega}|\boldsymbol{f}^y. \tag{5.38}$$

We write this as an algebraically equivalent variational problem on the same form as equation (5.4) by multiplying both sides with $\boldsymbol{y}'$, thus obtaining the discrete variational problem: Find $\boldsymbol{y} \in \boldsymbol{Y}$ such that

$$A_*^\delta(\boldsymbol{y}, \boldsymbol{y}') = \boldsymbol{f}^y(\boldsymbol{y}'), \tag{5.39}$$

where

$$A_*^\delta(\boldsymbol{y}, \boldsymbol{y}') = (\boldsymbol{A}_*\boldsymbol{y})\boldsymbol{y}' \qquad \text{and} \qquad \boldsymbol{f}^y(\boldsymbol{y}') = (|\boldsymbol{\omega}|\boldsymbol{f}^y, \boldsymbol{y}'). \tag{5.40}$$

Well-posedness of equation (5.39) again relies on continuity and coercivity of the bilinear forms, which will be shown below. However, we will need several technical tools, essentially mimicking discrete calculus rules, in the calculations. We summarize these relationships here:

**Lemma 5.9 [Calculus for FV-TPSA]**: The following relationships hold for the operators appearing in the FV-TPSA discretization:

1. Bounded averaging operator,
$$\left\|\tilde{\Xi}\boldsymbol{r}\right\|_{\bar{\mu}^{-1}}^{2} \leq \|\boldsymbol{r}\|_{\mu^{-1}}^{2}. \tag{5.41}$$
2. Orthogonal directional derivative decomposition,
$$\|\boldsymbol{\delta}^{-1}\Delta\boldsymbol{u}\|_{\bar{\mu}}^{2} = \|\boldsymbol{\delta}^{-1}\boldsymbol{n}\Delta\boldsymbol{u}\|_{\bar{\mu}}^{2} + \|\boldsymbol{\delta}^{-1}R^{\boldsymbol{n}}\Delta\boldsymbol{u}\|_{\bar{\mu}}^{2}. \tag{5.42}$$
3. The off-diagonal normal component terms are negative adjoints,
$$(\Delta|\boldsymbol{\varsigma}|\boldsymbol{n}\Xi\boldsymbol{u},\boldsymbol{p}) + \left(\Delta|\boldsymbol{\varsigma}|\boldsymbol{n}\tilde{\Xi}\boldsymbol{p},\boldsymbol{u}\right) = 0. \tag{5.43}$$
4. The off-diagonal rotation component terms are adjoints,
$$(\Delta|\boldsymbol{\varsigma}|R^{\boldsymbol{n}}\Xi\boldsymbol{u},\boldsymbol{r}) = \left(\Delta|\boldsymbol{\varsigma}|R^{\boldsymbol{n}}\tilde{\Xi}\boldsymbol{r},\boldsymbol{u}\right). \tag{5.44}$$

*Proof*: The above claims follow from algebraic manipulations of the stated operators, and are collected in Appendix A3. □

**Lemma 5.10 [Continuity of discrete elastic subsystem, $\vartheta = 0$]**: Subject to Assumptions 1.9, the bilinear form $A_{*}^{\delta}(\boldsymbol{y},\boldsymbol{y}')$ is continuous with respect to the norm in equation (5.37), and $\boldsymbol{f}^{y}(\boldsymbol{y}')$ is continuous with respect to the $L^2$-norm, with constants independent of $\delta$.

*Proof*: Since the problem is finite-dimensional, continuity is ensured. However, it remains to verify that the continuity constants are independent of $\delta$. From the definition of the $L^2$ norm (equation (A1.12)) and the Cauchy-Schwarz inequality, this holds for the linear term:

$$\boldsymbol{f}^{y}(\boldsymbol{y}') \leq \|\boldsymbol{f}^{y}\|\|\boldsymbol{y}'\|. \tag{5.45}$$

For the bilinear term the Cauchy-Schwarz inequality gives:

$$\begin{aligned} A_{*}^{\delta}(\boldsymbol{y},\boldsymbol{y}') \leq 2\|\boldsymbol{\delta}^{-1}\Delta^{*}\boldsymbol{u}\|_{\bar{\mu}}\|\boldsymbol{\delta}^{-1}\Delta^{*}\boldsymbol{u}'\|_{\bar{\mu}} + \left\|R^{\boldsymbol{n}}\tilde{\Xi}\boldsymbol{r}\right\|_{\bar{\mu}^{-1}}\|\boldsymbol{\delta}^{-1}\Delta^{*}\boldsymbol{u}'\|_{\bar{\mu}} + \|\boldsymbol{n}\tilde{\Xi}\boldsymbol{p}\|_{\bar{\mu}^{-1}}\|\boldsymbol{\delta}^{-1}\Delta^{*}\boldsymbol{u}'\|_{\bar{\mu}} + \\ \left\|R^{\boldsymbol{n}}\tilde{\Xi}\boldsymbol{r}'\right\|_{\bar{\mu}^{-1}}\|\boldsymbol{\delta}^{-1}\Delta^{*}\boldsymbol{u}\|_{\bar{\mu}} + \|\boldsymbol{n}\tilde{\Xi}\boldsymbol{p}'\|_{\bar{\mu}^{-1}}\|\boldsymbol{\delta}^{-1}\Delta^{*}\boldsymbol{u}\|_{\bar{\mu}} + \|\boldsymbol{r}\|_{\mu^{-1}}\|\boldsymbol{r}'\|_{\mu^{-1}} + \\ \|\Delta^{*}\boldsymbol{p}\|_{\boldsymbol{\delta}^{-1}\boldsymbol{\delta}^{\mu}}\|\Delta^{*}\boldsymbol{p}'\|_{\boldsymbol{\delta}^{-1}\boldsymbol{\delta}^{\mu}} + \|\boldsymbol{p}\|_{\lambda^{-1}}\|\boldsymbol{p}'\|_{\lambda^{-1}}. \end{aligned} \tag{5.46}$$

Here we used equations (5.43-5.44) of Lemma 5.9. We note that since $R^{\boldsymbol{n}}$ and $\Xi$ are bounded linear operators, and since $(\boldsymbol{\delta}^{-1}\boldsymbol{\delta}^{\mu})_k \lesssim \bar{\mu}_k^{-1}$ and $\delta_k^{\mu} = 0$ for boundary faces (given Dirichlet boundary conditions), all terms in (5.46) are bounded by terms in the norm given in equation (5.37), independent of $\delta \to 0$, thus the Lemma follows. □

In the analysis of the continuous problem, Section 5.1.1, we used the inf-sup properties of Stokes' equations, recalled in Lemma 5.3. A somewhat weaker result is available for the FV-TPSA, which we will now state and prove. This proof technique is motivated by the classic paper by Franca and Stenberg on stabilized methods [33], see also an earlier application to the analysis of finite volume methods for elasticity [28].

**Lemma 5.11 [Relaxed discrete Stokes' inf-sup]**: For any vector $\boldsymbol{p} \in \mathbb{R}^{|I_\mathcal{T}|}$, there exists a vector $\boldsymbol{u_p} \in \mathbb{R}^{n\times|I_\mathcal{T}|}$, such that:

$$\left(\Delta|\boldsymbol{\varsigma}|\boldsymbol{n}\Xi\boldsymbol{u_p},\boldsymbol{p}\right) \geq |\boldsymbol{p}|_{\mu^{-1}}\left(\beta_1|\boldsymbol{p}|_{\mu^{-1}} - \beta_2\|\Delta^*\boldsymbol{p}\|_{\boldsymbol{\delta}^{-1}\boldsymbol{\delta}^{\mu}}\right). \tag{5.47}$$

while at the same time:

$$\left(\Delta|\boldsymbol{\varsigma}|R^{\boldsymbol{n}}\Xi\boldsymbol{u_p},\boldsymbol{r}\right) \leq \beta_3|\boldsymbol{p}|_{\mu^{-1}}\|\boldsymbol{r}\|_{\mu^{-1}}. \tag{5.48}$$

The vector $\boldsymbol{u_p}$ further satisfies that:

$$\left\|\boldsymbol{\delta}^{-1}\Delta^*\boldsymbol{u_p}\right\|_{\overline{\mu}} \leq |\boldsymbol{p}|_{\mu^{-1}}. \tag{5.49}$$

*Proof*: The proof is by explicit construction. Let therefore $\boldsymbol{p}$ be given and denote by $p(x)$ the function that is constant on each cell $\omega_i \in \mathcal{T}$, taking the value of $p_i$. Then Lemma 5.3 gives the existence of a function $u_p \in H^1(\Omega,\mathbb{R}^n)$, with all the properties for the continuous Stokes' inf-sup. We now show that the cell-projection

$$\boldsymbol{u_p} = \pi_\omega u_p. \tag{5.50}$$

satisfies the properties stated in the Lemma. First note that equation (5.8) and the stability of the projection operator ensures (5.49), since:

$$|\boldsymbol{p}|_{\mu^{-1}} = |p(x)|_{\mu^{-1}} = \left\|\nabla u_p\right\|_{\mu} \geq \left\|\boldsymbol{\delta}^{-1}\Delta^*(\pi_\omega u_p)\right\|_{\overline{\mu}} = \left\|\boldsymbol{\delta}^{-1}\Delta^*\boldsymbol{u_p}\right\|_{\overline{\mu}}. \tag{5.51}$$

We now turn to equation (5.47).

$$\begin{aligned}\left(\Delta|\boldsymbol{\varsigma}|\boldsymbol{n}\Xi\boldsymbol{u_p},\boldsymbol{p}\right) &= \sum_{i\in I_\mathcal{T}} p_i \sum_{k\in\mathcal{N}_i}\int_{\varsigma_k}\Xi\boldsymbol{u_p}\cdot n_{k,i}\,dA \\ &= \textstyle\sum_{i\in I_\mathcal{T}} p_i \sum_{k\in N_i}\int_{\varsigma_k}\left(\Xi\boldsymbol{u_p} + u_p - \pi_{\varsigma,k}u_p\right)\cdot n_{k,i}\,dA.\end{aligned} \tag{5.52}$$

Where in the last term we have used that the $L^2$-projection onto the faces $\pi_{\varsigma,k}$ (defined in equation (A1.13)), preserves fluxes. Thus:

$$
\begin{aligned}
(\Delta|\boldsymbol{\varsigma}|\boldsymbol{n}\Xi\boldsymbol{u_p},\boldsymbol{p}) &= \sum_{i\in I_{\mathcal{T}}} p_i \sum_{k\in\mathcal{N}_i} \int_{\varsigma_k} \Xi_k \boldsymbol{u_p}\cdot n_{k,i}\,dA \\
&= (\nabla\cdot u_p, p) + \sum_{i\in I_{\mathcal{T}}} p_i \sum_{k\in\mathcal{N}_i} \int_{\varsigma_k} (\Xi_k\pi_\omega - \pi_{\varsigma,k})u_p\cdot n_{k,i}\,dA \\
&= (\nabla\cdot u_p, p) + \left(|\boldsymbol{\varsigma}|\boldsymbol{n}(\Xi\pi_\omega - \pi_\varsigma)u_p, \Delta^*\boldsymbol{p}\right) \\
&\geq \beta_1|\boldsymbol{p}|^2_{\mu^{-1}} - \left\|(\boldsymbol{\delta}\boldsymbol{\delta}^\mu)^{-\frac{1}{2}}\boldsymbol{n}(\Xi\pi_\omega - \pi_\varsigma)u_p\right\| \|\Delta^*\boldsymbol{p}\|_{\boldsymbol{\delta}^{-1}\boldsymbol{\delta}^\mu}.
\end{aligned}
\tag{5.53}
$$

Here $\beta_1$ is the inf-sup constant of the continuous case, see Lemma 5.3.

Furthermore, due to a self-similarity argument the difference between the mean of the cell interpolants and the face interpolants satisfies:

$$
\|(\boldsymbol{\delta}\boldsymbol{\delta}^\mu)^{-1/2}\boldsymbol{n}(\Xi\pi_\omega - \pi_\varsigma)u_p\| \leq \beta_2\|\nabla u_p\|_\mu = \beta_2|\boldsymbol{p}|_{\mu^{-1}}\,. \tag{5.54}
$$

where $\beta_2$ depends on the grid structure and the material constant $\mu$, but does not in general scale with $\delta$. Combining (5.52-5.54) gives (5.47).

Inequality (5.48) follows by a completely analogous calculation as equations (5.52-5.54), but exploiting the fact that $S\nabla u_p = 0$ to eliminate the first term in the resulting inequality, to obtain:

$$
(\Delta|\boldsymbol{\varsigma}|R^{\boldsymbol{n}}\Xi\boldsymbol{u_p},\boldsymbol{r}) \geq \beta_3|\boldsymbol{p}|_{\mu^{-1}}\|\Delta^*\boldsymbol{r}\|_{\boldsymbol{\delta}^{-1}\boldsymbol{\delta}^\mu}. \tag{5.55}
$$

Inequality (5.48) is now finally obtained by the inverse inequality [23]:

$$
\|\Delta^*\boldsymbol{r}\|_{\boldsymbol{\delta}^{-1}\boldsymbol{\delta}^\mu} \leq \|\boldsymbol{r}\|_{\mu^{-1}}. \tag{5.56}
$$

□

We are now ready to prove the main stability result for FV-TPSA.

**Theorem 5.12 [Coercivity of discrete elastic subsystem, $\vartheta = 0$]**: Subject to Assumption 1.9 the bilinear form $A_*^\delta(\boldsymbol{y},\boldsymbol{y}')$ satisfies:

$$
\inf_{\boldsymbol{y}\in\boldsymbol{Y}}\sup_{\boldsymbol{y}'\in\boldsymbol{Y}}\frac{A_*^\delta(\boldsymbol{y},\boldsymbol{y}')}{\|\boldsymbol{y}\|_*\|\boldsymbol{y}'\|_*} \gtrsim 1.
$$

where the hidden constant is robust in the sense of Definition 1.6, and independent of $\delta\to 0$.

*Proof*: The proof follows the structure of Lemma 5.4, but is adapted to the available discrete calculus identities (Lemma 5.9) and the relaxed Stokes' inf-sup condition (Lemma 5.11). Thus, as before, we show the condition by an explicit construction. To this end, let $\boldsymbol{y} = (\boldsymbol{u},\boldsymbol{r},\boldsymbol{p}) \in \boldsymbol{Y}$ be given, and choose $\boldsymbol{y}' = (\boldsymbol{u} + \alpha\boldsymbol{u_p},\boldsymbol{r},\boldsymbol{p})$, where $\alpha$ is a free parameter and $\boldsymbol{u_p}$ depends on $\boldsymbol{p}$ according to Lemma 5.11. In view of Lemma 5.9, we then calculate

$$A_*^\delta(\boldsymbol{y},\boldsymbol{y}') = 2\|\boldsymbol{\delta}^{-1}\Delta^*\boldsymbol{u}\|_{\bar{\mu}}^2 + 2\big(\Delta|\boldsymbol{\varsigma}|R^{\boldsymbol{n}}\tilde{\Xi}\boldsymbol{r},\boldsymbol{u}\big) + \|\boldsymbol{r}\|_{\mu^{-1}}^2 + \|\Delta^*\boldsymbol{p}\|_{\boldsymbol{\delta}^{-1}\boldsymbol{\delta}^{\mu}}^2 + \|\boldsymbol{p}\|_{\lambda^{-1}}^2 + \alpha\Big(2\big(|\boldsymbol{\varsigma}|\boldsymbol{\delta}^{-1}\bar{\boldsymbol{\mu}}\Delta^*\boldsymbol{u_p},\Delta^*\boldsymbol{u}\big) + \big(\Delta|\boldsymbol{\varsigma}|R^{\boldsymbol{n}}\Xi\boldsymbol{u_p},\boldsymbol{r}\big) + \big(\Delta|\boldsymbol{\varsigma}|\boldsymbol{n}\Xi\boldsymbol{u_p},\boldsymbol{p}\big)\Big). \tag{5.57}$$

Applying now equations (5.47-5.49) from Lemma 5.11, we obtain for $0 < \gamma_1 < 2$:

$$\begin{aligned} A_*^\delta(\boldsymbol{y},\boldsymbol{y}') \geq{}& (2-\gamma_1)\|\boldsymbol{\delta}^{-1}\Delta^*\boldsymbol{u}\|_{\bar{\mu}}^2 + \gamma_1\big(\|\boldsymbol{\delta}^{-1}\boldsymbol{n}\Delta^*\boldsymbol{u}\|_{\bar{\mu}}^2 + \|\boldsymbol{\delta}^{-1}R^{\boldsymbol{n}}\Delta^*\boldsymbol{u}\|_{\bar{\mu}}^2\big) \\ &- 2\|\boldsymbol{\delta}^{-1}R^{\boldsymbol{n}}\Delta^*\boldsymbol{u}\|_{\bar{\mu}}\big\|\tilde{\Xi}\boldsymbol{r}\big\|_{\bar{\mu}^{-1}} + \|\boldsymbol{r}\|_{\mu^{-1}}^2 + \|\boldsymbol{p}\|_{\lambda^{-1}}^2 + \|\Delta^*\boldsymbol{p}\|_{\boldsymbol{\delta}^{-1}\boldsymbol{\delta}^{\mu}}^2 \\ &+\alpha\Big(-2\|\boldsymbol{\delta}^{-1}\Delta^*\boldsymbol{u}\|_{\bar{\mu}}|\boldsymbol{p}|_{\mu^{-1}} - \beta_3|\boldsymbol{p}|_{\mu^{-1}}\|\boldsymbol{r}\|_{\mu^{-1}} + |\boldsymbol{p}|_{\mu^{-1}}\big(\beta_1|\boldsymbol{p}|_{\mu^{-1}} - \beta_2\|\Delta^*\boldsymbol{p}\|_{\boldsymbol{\delta}^{-1}\boldsymbol{\delta}^{\mu}}\big)\Big). \end{aligned} \tag{5.58}$$

Furthermore, we note that using equation (5.41) from Lemma 5.9, together with Young's inequality with weight $\gamma_2 > 0$, we obtain:

$$2\|\boldsymbol{\delta}^{-1}R^{\boldsymbol{n}}\Delta^*\boldsymbol{u}\|_{\bar{\mu}}\big\|\tilde{\Xi}\boldsymbol{r}\big\|_{\bar{\mu}^{-1}} \leq \gamma_2\|\boldsymbol{\delta}^{-1}R^{\boldsymbol{n}}\Delta^*\boldsymbol{u}\|_{\bar{\mu}}^2 + \frac{1}{\gamma_2}\big\|\tilde{\Xi}\boldsymbol{r}\big\|_{\bar{\mu}^{-1}}^2 \leq \gamma_2\|\boldsymbol{\delta}^{-1}R^{\boldsymbol{n}}\Delta^*\boldsymbol{u}\|_{\bar{\mu}}^2 + \frac{1}{\gamma_2}\|\boldsymbol{r}\|_{\mu^{-1}}^2. \tag{5.59}$$

Applying Young's inequality also to the remaining terms containing products of norms, we then obtain with positive weights $\gamma_3$, $\gamma_4$ and $\gamma_5$:

$$A_*^\delta(\boldsymbol{y},\boldsymbol{y}') \geq (2-\gamma_1-\alpha\gamma_3)\|\boldsymbol{\delta}^{-1}\Delta^*\boldsymbol{u}\|_{\bar{\mu}}^2 + (\gamma_1-\gamma_2)\|\boldsymbol{\delta}^{-1}R^{\boldsymbol{n}}\Delta^*\boldsymbol{u}\|_{\bar{\mu}}^2 + \Big(1-\frac{\alpha\beta_3}{2\gamma_4}-\frac{1}{\gamma_2}\Big)\|\boldsymbol{r}\|_{\mu^{-1}}^2 + \Big(1-\frac{\alpha\beta_2}{2\gamma_5}\Big)\|\Delta^*\boldsymbol{p}\|_{\boldsymbol{\delta}^{-1}\boldsymbol{\delta}^{\mu}}^2 + \|\boldsymbol{p}\|_{\lambda^{-1}}^2 + \bigg(\alpha\beta_1 - \alpha\Big(\frac{1}{\gamma_3}+\frac{\gamma_4\beta_3}{2}+\frac{\gamma_5\beta_2}{2}\Big)\bigg)|\boldsymbol{p}|_{\mu^{-1}}^2. \tag{5.60}$$

We eliminate non-essential terms, and in general simplify, by setting

$$\gamma_5 = \frac{\alpha\beta_2}{2} \qquad \text{and} \qquad \gamma_2 = \gamma_1. \tag{5.61}$$

Then:

$$A_*^\delta(\boldsymbol{y},\boldsymbol{y}') \geq (2-\gamma_1-\alpha\gamma_3)\|\boldsymbol{\delta}^{-1}\Delta^*\boldsymbol{u}\|_{\bar{\mu}}^2 + \Big(1-\frac{\alpha\beta_3}{2\gamma_4}-\frac{1}{\gamma_1}\Big)\|\boldsymbol{r}\|_{\mu^{-1}}^2 + \|\boldsymbol{p}\|_{\lambda^{-1}}^2 + \alpha\Big(\beta_1-\frac{1}{\gamma_3}-\frac{\gamma_4\beta_3}{2}-\frac{\alpha\beta_2^2}{4}\Big)|\boldsymbol{p}|_{\mu^{-1}}^2 \tag{5.62}$$

Now by setting e.g.

$$\gamma_1 = \frac{3}{2}, \qquad \gamma_4 = \sqrt{\alpha}\,, \qquad \gamma_3 = \sqrt{\alpha^{-1}}, \tag{5.63}$$

we obtain:

$$A_*^\delta(\boldsymbol{y},\boldsymbol{y}') \geq \Big(\frac{1}{2}-\sqrt{\alpha}\Big)\|\boldsymbol{\delta}^{-1}\Delta^*\boldsymbol{u}\|_{\bar{\mu}}^2 + \Big(\frac{1}{3}-\frac{\sqrt{\alpha}\beta_3}{2}\Big)\|\boldsymbol{r}\|_{\mu^{-1}}^2 + \|\boldsymbol{p}\|_{\lambda^{-1}}^2 + \alpha\Big(\beta_1-\sqrt{\alpha}-\frac{\sqrt{\alpha}\beta_3}{2}-\frac{\alpha\beta_2^2}{4}\Big)|\boldsymbol{p}|_{\mu^{-1}}^2. \tag{5.64}$$

We can now finally choose $\alpha$ sufficiently small, depending only on the three $\beta_i$, which are robust for non-degenerate $\mu$, as required by Assumption 1.9. It is then clear that all parentheses in (5.64) are positive. As all constants appearing in (5.64) are robust, so is the choice of $\alpha$. Following an application of the discrete Poincaré inequality (equation (A1.15)), and the fact that $\alpha$ is bounded and thus $\|\boldsymbol{y}'\|_\circ \lesssim \|\boldsymbol{y}\|_\circ$, concludes the proof. □

**Corollary 5.13 [Well-posedness of the FV-TPSA discretization for the elastic subsystem]**: Subject to Assumptions 1.9, equation (5.39) is well-posed, with solution $\boldsymbol{y} \in \boldsymbol{Y}$ satisfying

$$\|\boldsymbol{y}\|_* \lesssim \|\boldsymbol{f}^y\|$$

where the constants are robust in the sense of Definition 1.6, and do not depend on $\delta \to 0$.

*Proof*: The corollary follows Lemma 5.10, Theorem 5.12 and standard saddle-point theory [3] □.

## 5.2.2 Full poromechanical system

Having developed the well-posedness for the elastic subsystem, the well-posedness of the full system follows from an application of Lemma 5.7. Following the notation of Section 5.2.1, we consider the FV-TPSA discretization stated in equation (4.4). By designating the discretization matrix of the poromechanical system by $\boldsymbol{A}_\circ$, the FV-TPSA discretization can be written compactly as:

$$\boldsymbol{A}_\circ \boldsymbol{z} = |\boldsymbol{\omega}| \boldsymbol{f}^z. \tag{5.65}$$

We again write this as an algebraically equivalent variational problem by multiplying both sides with $\boldsymbol{z}'$, thus obtaining the discrete variational problem: Find $\boldsymbol{z} \in \boldsymbol{Z}$ such that

$$A_\circ^\delta(\boldsymbol{z}, \boldsymbol{z}') = \boldsymbol{f}^z(\boldsymbol{z}'), \tag{5.66}$$

where

$$A_\circ^\delta(\boldsymbol{z}, \boldsymbol{z}') = (\boldsymbol{A}_\circ \boldsymbol{z})\boldsymbol{z}' \qquad \text{and} \qquad \boldsymbol{f}^z(\boldsymbol{z}') = (|\boldsymbol{\omega}|\boldsymbol{f}^z, \boldsymbol{z}'). \tag{5.67}$$

**Lemma 5.14 [Continuity of discrete poromechanical system]**: Subject to Assumption 1.11, the bilinear form $A_\circ^\delta(\boldsymbol{z}, \boldsymbol{z}')$ is continuous with respect to the norm in equation (5.37), and $\boldsymbol{f}^z(\boldsymbol{z}')$ is continuous with respect to the $L^2$-norm with constants independent of $\delta$.

*Proof*: The proof is a straight-forward application of the definition of $A_\circ^\delta$, the Cauchy-Schwartz inequality, and Lemma 5.9. □

**Lemma 5.15 [Coercivity of discrete coupled poroelasticity]**: Subject to Assumption 1.11, the bilinear form $A_\circ^\delta(\boldsymbol{z}, \boldsymbol{z}')$ satisfies:

$$\inf_{\boldsymbol{z}\in\boldsymbol{Z}} \sup_{\boldsymbol{z}'\in\boldsymbol{Z}} \frac{A_\circ^\delta(\boldsymbol{z},\boldsymbol{z}')}{\|\boldsymbol{z}\|_*\|\boldsymbol{z}'\|_*} \gtrsim 1. \tag{5.68}$$

where the hidden constant is robust in the sense of Definition 1.6, and independent of $\delta \to 0$.

*Proof*: We apply Lemma 5.7. From the definition of $A_\circ^\delta$, it is clear that it is on the form of (5.25), and that equation (5.28) holds. Moreover, condition 1) of the proof was shown in equation (5.64) of the proof of Theorem 5.12. It remains to show that condition 2) holds, e.g. that the Darcy system satisfies

$$\inf_{\boldsymbol{w}\in\boldsymbol{W}} \frac{C^\delta(\boldsymbol{w},\boldsymbol{w}) - \|\boldsymbol{w}\|_\eta^2}{\|\delta^{-1}\bar{\boldsymbol{\kappa}}\boldsymbol{w}\|_\kappa^2} \gtrsim 1, \tag{5.69}$$

where from equation (4.4) we identify

$$C^\delta(\boldsymbol{w}, \boldsymbol{w}') = \left(\Delta|\boldsymbol{\varsigma}|\boldsymbol{\delta}^{-1}\bar{\boldsymbol{\kappa}}\Delta^* \boldsymbol{w}, \boldsymbol{w}'\right) + (|\boldsymbol{\omega}|\boldsymbol{\eta}\boldsymbol{w}, \boldsymbol{w}'). \tag{5.70}$$

However, for any $\boldsymbol{w}$, equation (5.69) holds with equality, thus all conditions of Lemma 5.7 hold, and inequality (5.68) is thus satisfied. □

**Theorem 5.16 [Well-posedness of the FV-TPSA]**: Subject to Assumption 1.11, equation (5.66) is well-posed with solution $\boldsymbol{z} \in \boldsymbol{Z}$ satisfying

$$\|\boldsymbol{z}\|_\circ \lesssim \|\boldsymbol{f}^z\|,$$

where the constants are robust in the sense of Definition 1.6, and do not depend on $\delta \to 0$.

*Proof*: The theorem follows Lemma 5.14 and 5.15 and standard saddle-point theory [3] □.

## 5.3 Consistency and convergence of FV-TPSA

It is clear that two-point approximations cannot provide a consistent approximation to a normal derivative if the vector obtained by subtracting the coordinate vector of the two points is not parallel to the normal vector. In the present setting, this appears in Equation (A2.3), which based on a Taylor series expansion can be seen to have an approximation error of

$$n_{i,j}^T \nabla u|_{x_{\{i,j\}}} = \frac{u|_{x_{\{i,j\}}} - u_i}{(x_{\{i,j\}} - x_i) \cdot n_{i,j}} + \mathcal{O}(1) \left| n_{i,j} \times \nabla u|_{x_{\{i,j\}}} \right| + \mathcal{O}(\delta) \left( n_{i,j}^T \nabla \left( n_{i,j}^T \nabla u \right) \Big|_{x_{\{i,j\}}} \right) + \mathcal{O}(\delta^2). \tag{5.71}$$

As such, no two-point scheme will provide a consistent numerical flux, unless the grid is face-orthogonal, in the sense of Definition 2.3.

It is therefore *a priori* clear that there will exist classes of grids (such as e.g. parallelograms regular tiling of non-equilateral triangles) for which the FV-TPSA method is not consistent, just like its scalar counterpart TPFA (for an in-depth discussion, see [29]). Consistency must therefore be established on a smaller class of grids. *As a consequence, for this section, we only consider face-orthogonal grids*.

The consistency and convergence of two-point approximation schemes for elliptic PDEs is carefully treated in the classic work of Eymard, Gallouët and Herbin [23]. In particular, much of the results of Section 3.1.4 to 3.1.6 of that work directly applies to the currently proposed methods. In the interest of space, we will not reproduce their arguments as adapted to the current context, but state the main results in the following theorem.

**Theorem 5.17 [Convergence of the FV-TPSA]**: For admissible face-orthogonal grids, and spatially constant material parameters, the solution $z$ of equation (5.22) and the solution $\boldsymbol{z}_\delta$ of equation (5.66) (for a given grid $\mathcal{T}_\delta$ in a grid sequence indexed by $\delta$), satisfy whenever $z$ is sufficiently smooth:

$$\|\boldsymbol{z}_\delta - z(\boldsymbol{x}_\omega)\|_\circ \lesssim \|f^z\|\delta. \tag{5.72}$$

where $\boldsymbol{x}_\omega$ is the vector of cell centers, and the hidden constant is robust in the sense of Definition 1.6, and does not depend on $\delta \to 0$.

*Proof*. Existence of $\boldsymbol{z}_\delta$ and $z$ is guaranteed by Theorems 5.8 and 5.16. The closeness claimed in the Theorem follows the same arguments as in the proof of Theorem 3.3 in [23]. □

**Corollary 5.18 [Convergence of the FV-TPSA stress and flux]**: Whenever Theorem 5.17 applies, the numerical stresses and flux $(\boldsymbol{\sigma}_\delta, \boldsymbol{\tau}_\delta, \boldsymbol{\chi}_\delta)$, defined from $\boldsymbol{z}_\delta$ by equation (4.3), also converge to the continuous stresses and flux $(\sigma, \tau, \chi)$, defined from $z$ by equation (1.15), satisfy

$$\left\|\boldsymbol{\sigma}_\delta - \sigma(\boldsymbol{x}_\varsigma)\right\|_{\overline{\mu}^{-1}}^2 + \left\|\boldsymbol{\chi}_\delta - \chi(\boldsymbol{x}_\varsigma)\right\|_{\overline{\kappa}^{-1}}^2 \lesssim \|f^z\|\delta. \tag{5.73}$$

where $\boldsymbol{x}_\varsigma$ is the vector of cell centers, and the hidden constant is robust in the sense of Definition 1.6, and does not depend on $\delta \to 0$.

*Proof*. From equation (4.3), the stated stresses and flux are algebraically related to quantities appearing in the norm of Theorem 5.17. Keeping in mind that by Assumption

1.11, $\mu$ is bounded from below all constants are reflected in the norms stated in the Corollary. □

# 6. Numerical verification

We verify the performance of the finite volume TPSA discretization through a series of numerical experiments, with emphasis on probing the robustness of the method in the sense of Definition 1.6.

From the perspective of grids, we identify four grids that illustrate the performance and robustness of the method, we refer to these as Grid Types (GT) 1 through 4, see Figure 6.1 for an illustration:

GT1: *Symmetric grids*: For grids with high degree of local symmetry, second-order convergence is often seen for discretizations of elliptic-type partial differential equations, despite only first-order convergence being proved (see e.g. [34, 35]). We explore this with a regular Cartesian grid sequence.

GT2: *Unstructured face-orthogonal grids*: For face-orthogonal grids without additional symmetry, Theorem 5.16 should apply. We explore this with a sequence of prismatic extensions of triangular grids provided by Gmsh [36].

GT3: *Asymptotically face-orthogonal grids*: Theorem 5.16 requires face-orthogonal grid sequences, but it is natural to expect good results also for minor deviations from face-orthogonality. We explore this with a prismatic extrusion of a 2D logically Cartesian grid sequence obtained by a $\delta^{3/2}$ perturbation of the internal nodes.

GT4: *Non-face-orthognal grids:* We do not expect convergence on very bad grids, however, it is reasonable to still expect a stable approximate solution, in view of Theorem 5.16. We explore this with a prismatic extrusion of a 2D logically Cartesian grid sequence obtained by a $\delta$ perturbation of the internal nodes.

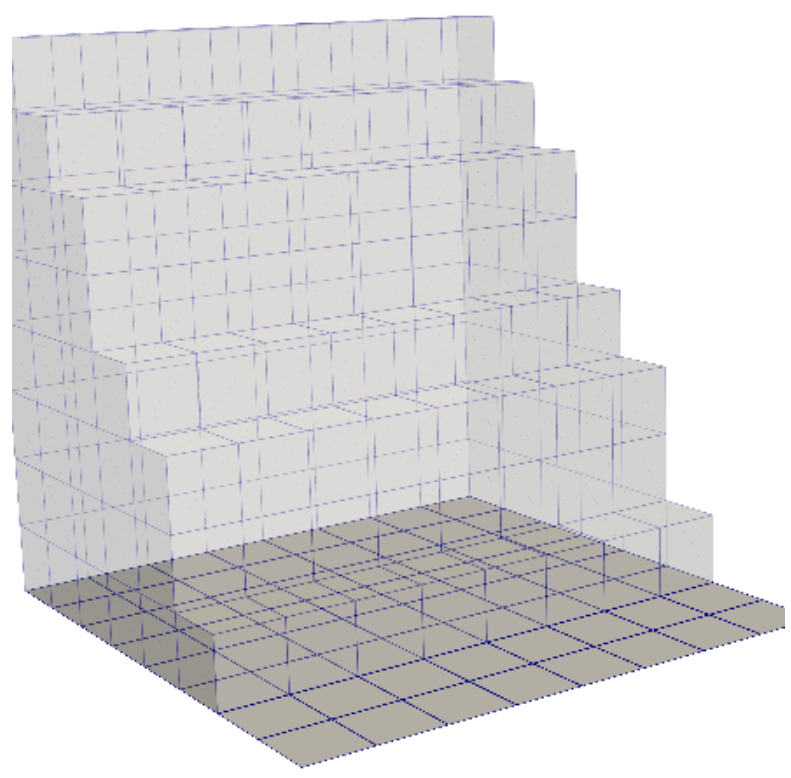

GT1

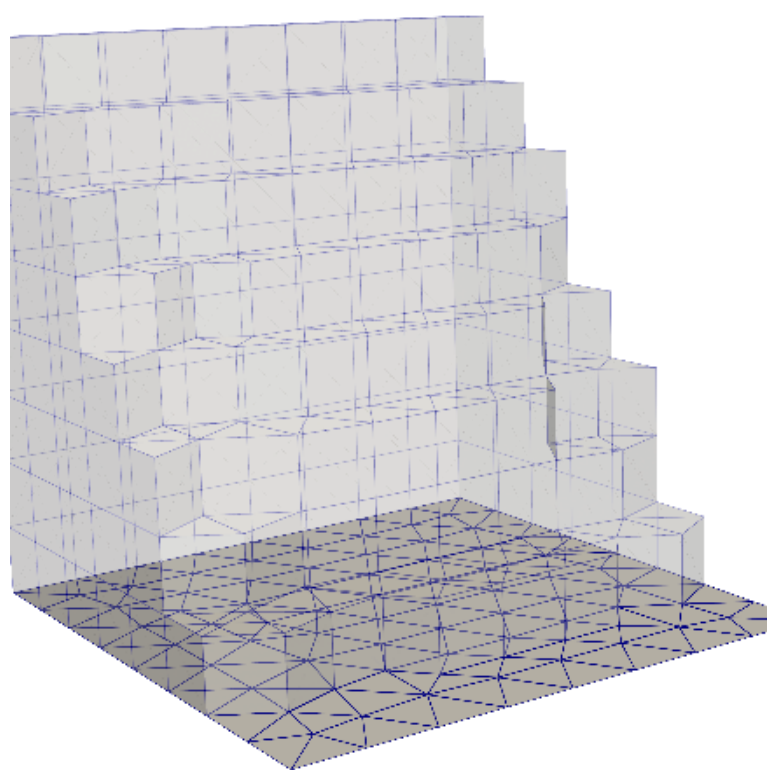

GT2

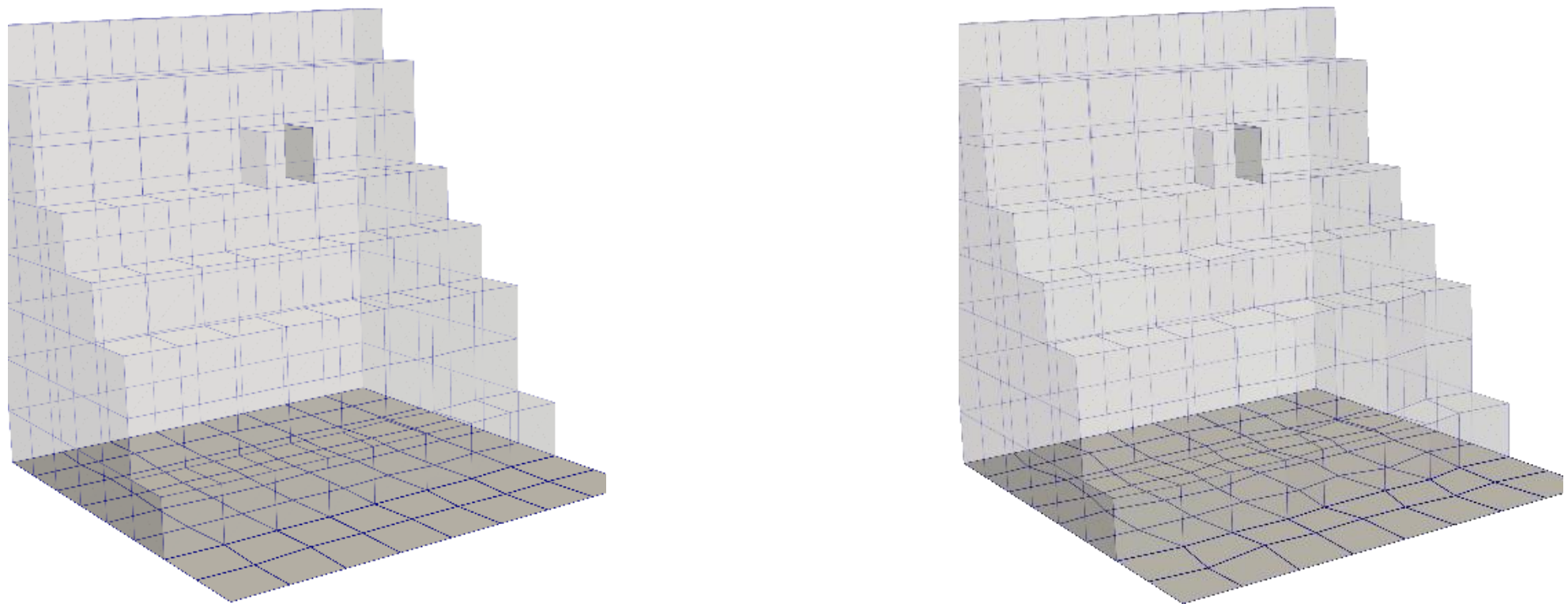

**Figure 6.1:** Illustration of the grid types used in the convergence test. The grids have been cut to better show the perturbations. Upper left: Unperturbed Cartesian grid (GT1). Upper right: Simplex grid (GT2). Lower left: Cartesian grid with $\delta^{3/2}$ perturbations (GT3). Lower right: Cartesian grid with $\delta$ perturbations (GT4).

We note that the GT4 is a quite challenging grid type in general, and also more advanced methods such as the symmetric variants of the so-called multipoint mixed-finite element methods lose convergence on these grids [37, 12]. In our experience, GT4 is also indicative of the performance on general 3D simplicial grids, whenever these are sufficiently deformed from an isosceles so that the circumcenter does not lie within the cell.

All experiments are on the unit cube in 3D. The TPSA and TPFA discretizations are implemented in the open-source simulation code PorePy [38]. The resulting linear systems are solved with GMRES preconditioned with smoothed aggregation algebraic multigrid methods provided by PyAMG [39]. Right hand side terms are computed by a first order quadrature rule.

Following Theorem 5.71, we compute the error of the elastic system, referred to as $e_*$, by the norm defined in equation (5.36), while the error for the poromechanical system, $e_\circ$ is computed according to (5.37). The reported errors are scaled with the norm of the manufactured solution on the finest grid level.

Due to space constraints, we do not include comparison to existing discretizations such as the mixed finite element method and the multi-point stress finite volume method. A 2D comparison to these methods has been communicated separately [40].

The numerical experiments can be reproduced using a Docker image available at [41], using runscripts located at [42].

## 6.1 Linearized elasticity and Stokes

We first consider the case of a purely elastic material ($\vartheta = 0$) with emphasis on robustness in the incompressible limit, and therefore set $\lambda = \{1, 10^2, 10^4, 10^{10}\}$. We report the convergence towards a manufactured expression, slightly modified from [14], and chosen to conform to zero Dirichlet boundary conditions:

$$u = \sum_{i=1}^{3} e_i \left( \frac{\partial \psi}{\partial x_{i+1}} - \frac{\partial \psi}{\partial x_{i-1}} \right), \qquad \psi = \prod_{i=1}^{3} \sin^2(\pi x_i),$$

$$r = \sum_{i=1}^{3} 100 x_i (1 - x_i) \sin(\pi x_{i+1}) \sin(\pi x_{i-1})\, e_i \;, \quad p = 0. \tag{6.1}$$

Here, the indices are taken as modulo 3, and $e_i$ is a unit vector along the $i$-th coordinate axis.

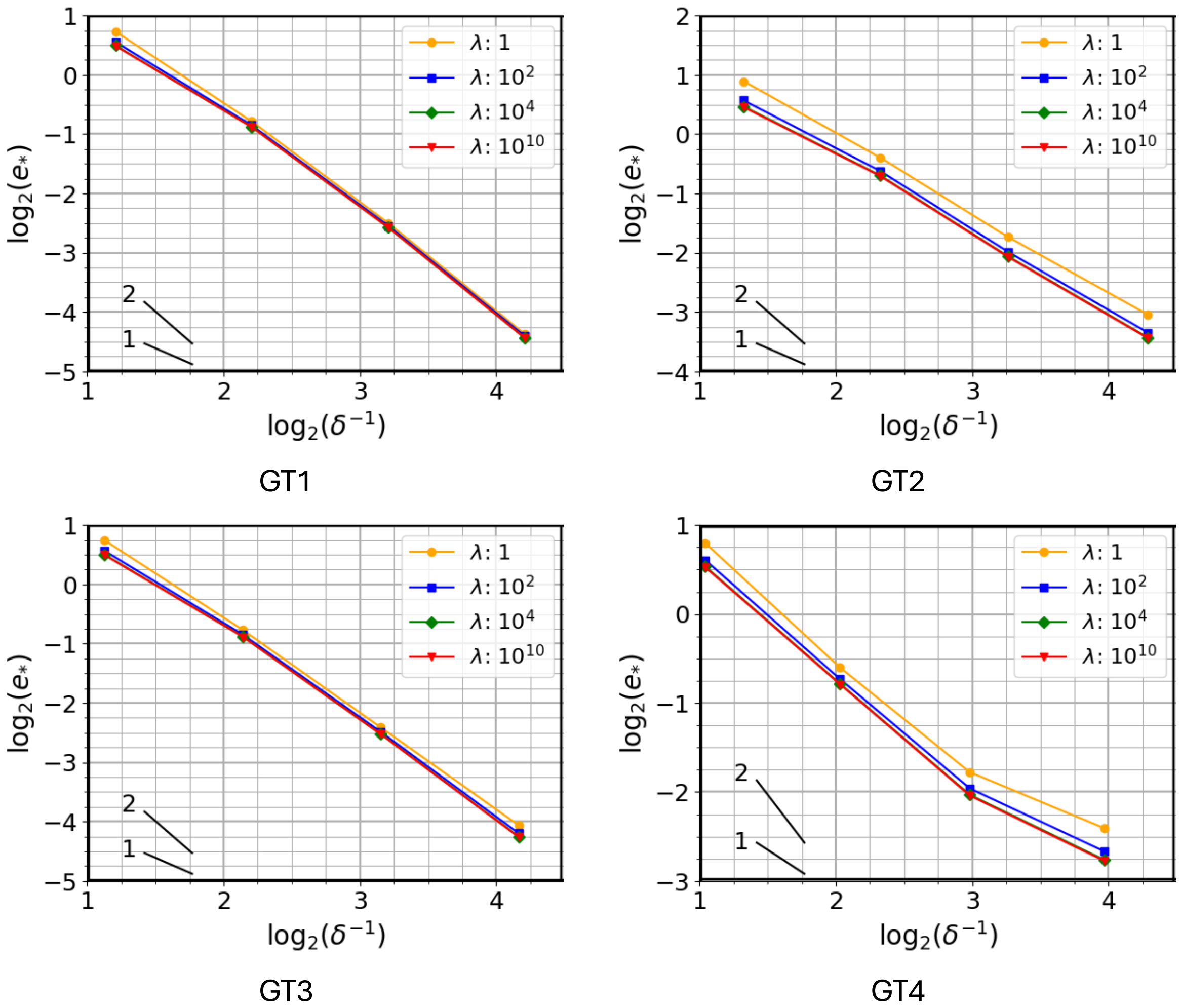

**Figure 6.2:** Convergence for an elastic material, considering robustness in the incompressible limit. Upper left: Unperturbed Cartesian grid (GT1). Upper right: Prismatically extended simplex grid (GT2). Lower left: Cartesian grid with $\delta^{3/2}$ perturbations (GT3). Lower right: Cartesian grid with $\delta$ perturbations (GT4).

The errors under grid refinement are shown in Figure 6.2 for grid types GT1-GT4. All results are consistent with the theoretical expectations: The discretization is stable for all grids, independently of $\lambda$, and converges for grids GT1, GT2 and GT3. Moreover, second order convergence is observed for GT1, first order for GT2 and order 1.5 for GT3. For GT4, the non-face-orthogonal grid, the solution is stable, but the observed convergence deteriorates under grid refinement and further loss of convergence can be expected for finer grids. As can be seen from Figure 6.3, the displacement, which in a sense is the primary variable for this problem, appears to have improved convergence rates and exhibits second order convergence on GT1-GT3, while also being more accurately approximated on GT4.

All results are fully robust for the parameter range considered. In 2D computations and applying direct solvers, we have also considered the case of $\lambda^{-1} = 0$, and verified that the solution is convergent, appreciating that the solid pressure is unique only up to an additive constant, as reflected in the norm (5.1). We expect this to hold also in 3D, however, our preconditioner is not implemented to deal with the undetermined solid pressure, and we therefore limit the investigation to $\lambda = 10^{10}$.

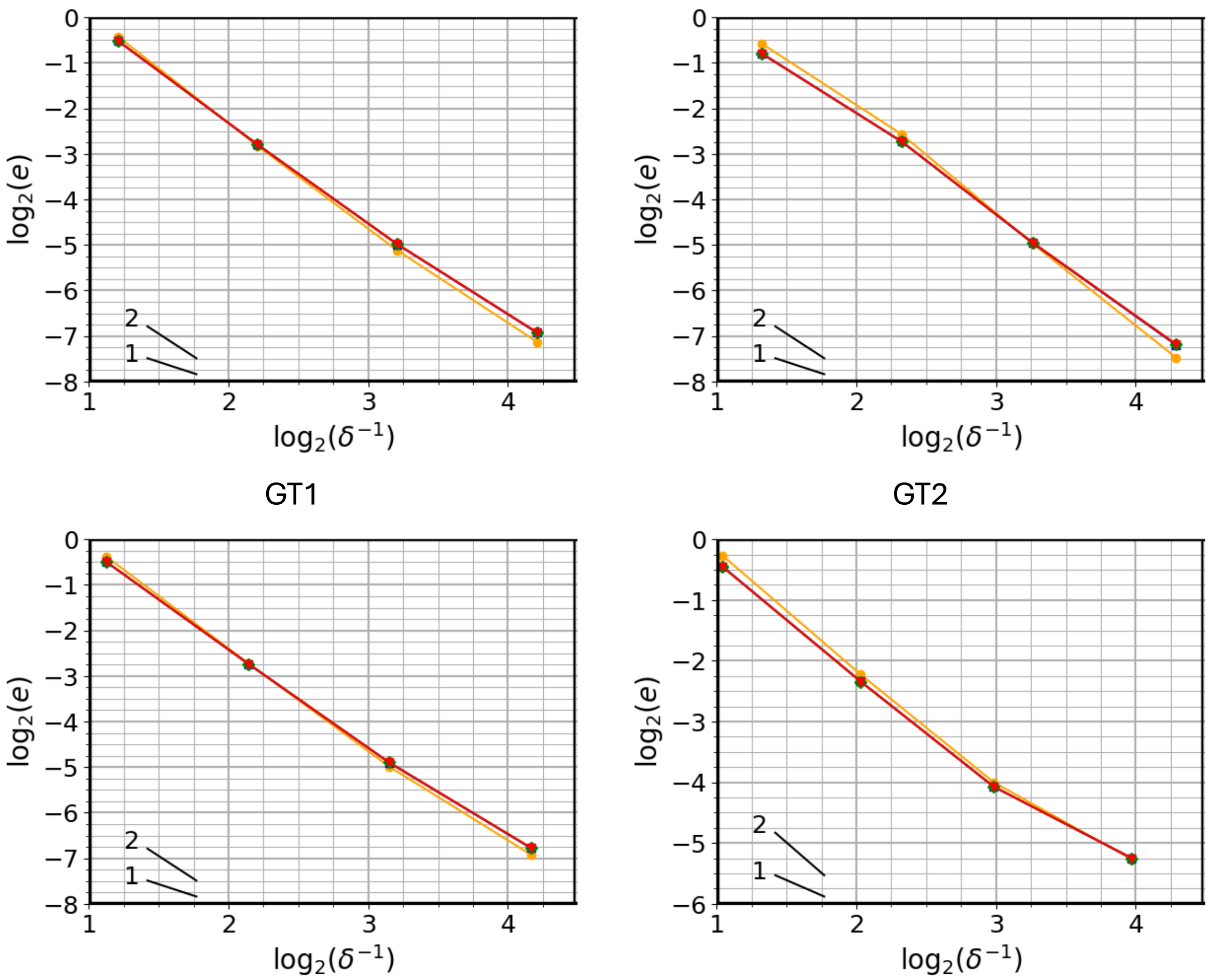

GT3 GT4

**Figure 6.3:** Convergence of the displacement for an elastic material, considering robustness in the incompressible limit. The color coding is the same as in Figure 6.2. Upper left: Unperturbed Cartesian grid (GT1). Upper right: Prismatically extended simplex grid (GT2). Lower left: Cartesian grid with $\delta^{3/2}$ perturbations (GT3). Lower right: Cartesian grid with $\delta$ perturbations (GT4).

## 6.2 Material heterogeneity

To test robustness in the presence of material heterogeneities, we consider a unit cube domain with a material heterogeneity which can be described as follows: We define a characteristic function

$$\chi = \begin{cases} 1 & \min(x_1, x_2, x_3) > 1/2 \\ 0 & \text{otherwise} \end{cases}, \tag{6.2}$$

and define spatially heterogeneous Lame parameters $\mu = (1-\chi) + \chi\beta$ and $\lambda = \mu$. Whenever $\beta \neq 1$, the material parameters are discontinuous, and we consider the cases $\beta = \{10^{-4}, 1, 10^{4}\}$. The analytical solution is set to

$$u = \frac{\sum_{i=1}^{3} \sin(2\pi x_i)\, x_{i+(}(1-x_{i+1})x_{i-1}(1-x_{i-1})e_i}{(1-\chi)+\chi\beta}$$

$$r = (e_1 + e_2 + e_3)\prod_{i=1}^{3} x_i\left(x_i - \frac{1}{2}\right)(x_i - 1)$$

$$p = \sin(2\pi x_1)\, x_2(1-x_2)x_3(1-x_3)$$

(6.3)

The solution is constructed so that the displacement gradient has a sharp discontinuity at the material contrast, but with continuous displacement and normal component of the stress. For all grid types, the mesh is constructed to resolve the material discontinuity.

The errors under grid refinement are shown in Figure 6.4, and are analogous to those reported in Section 6.1. Again, the discretization is stable for all grids, independently of the contrast in material parameters. In general, we observe second order convergence on GT1 and first order on GT2 and GT3. For GT2, the refinement levels obtained from Gmsh contain some arbitrariness, and the grid quality is therefore not as uniform as for the other grid types. This is in particular reflected for $\beta = 10^{-4}$, where the convergence is below one for the first two refinement steps, followed by rapid convergence at the final refinement. For GT4 the error does not decrease under grid refinement but remains stable for all values of $\beta$. As can be seen in Figure 6.5, the displacement variable appears to exhibit second order convergence for GT1 and above one for GT3 and GT4. For GT2 we in general observe first order convergence, with some irregularity as discussed above.

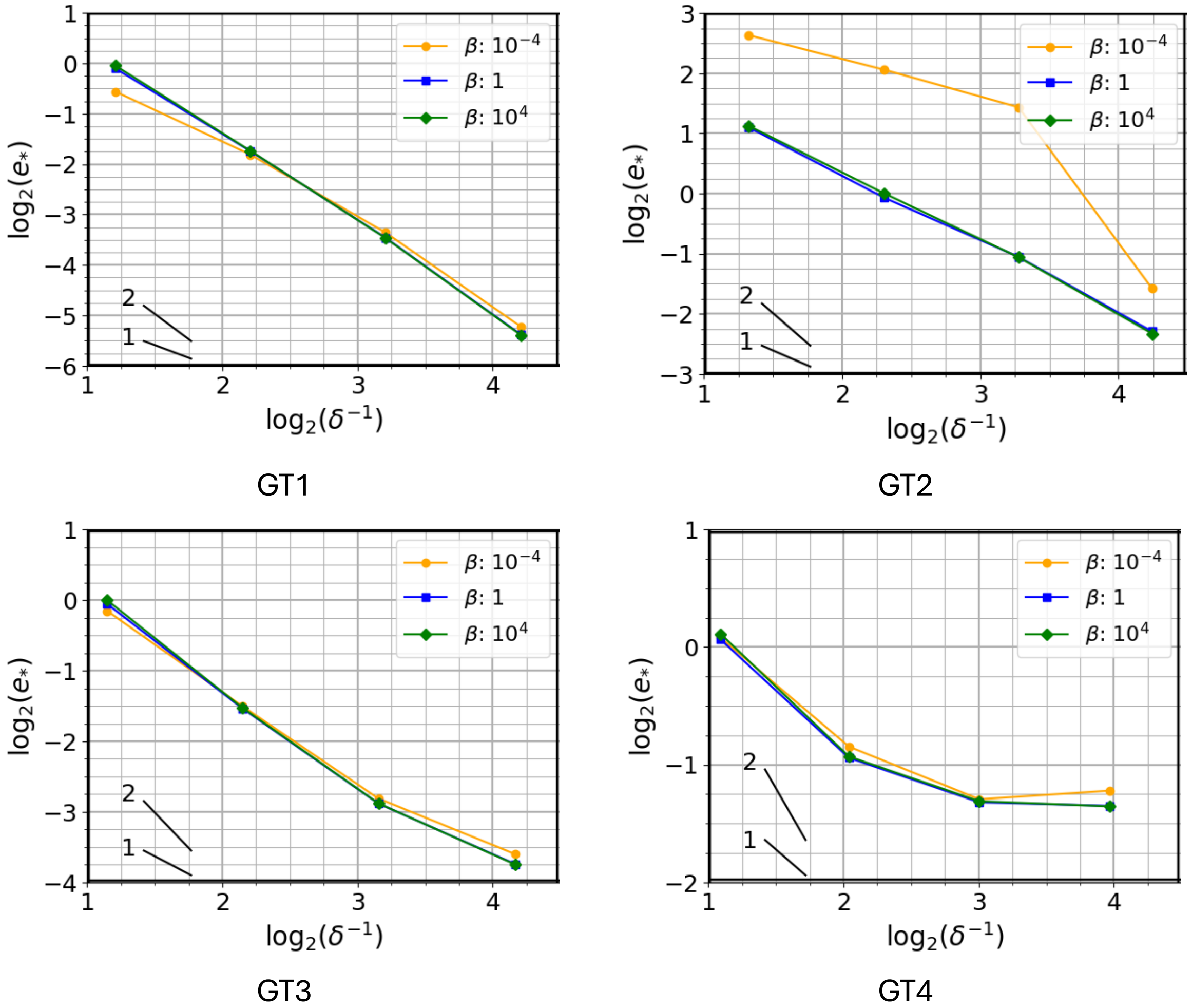

**Figure 6.4:** Convergence for an elastic material, considering robustness under a material heterogeneity. Upper left: Unperturbed Cartesian grid (GT1). Upper right: Prismatically extended simplex grid (GT2). Lower left: Cartesian grid with $\delta^{3/2}$ perturbations (GT3). Lower right: Cartesian grid with $\delta$ perturbations (GT4).

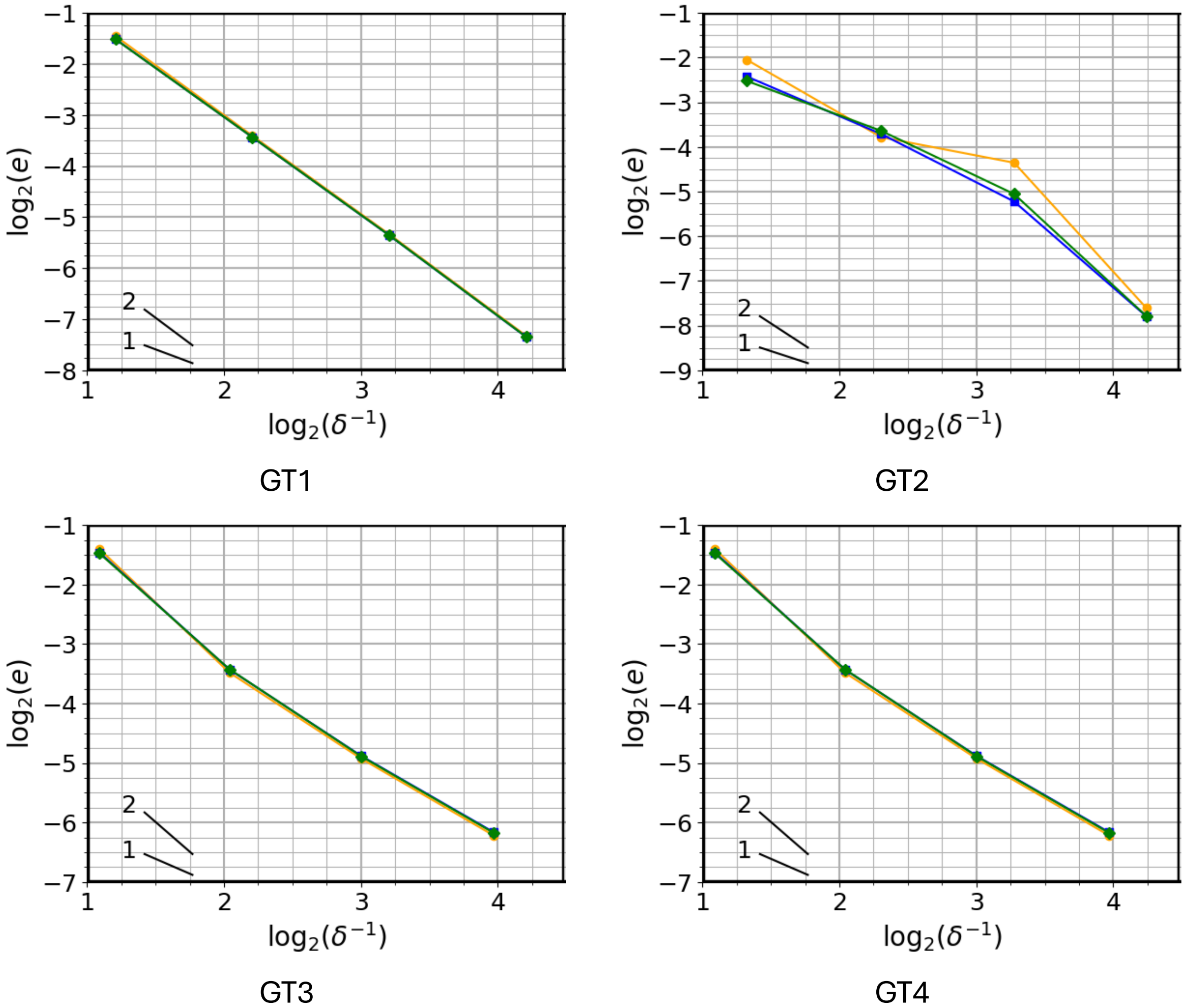

**Figure 6.5:** Convergence for an elastic material, considering robustness under a material heterogeneity. The color coding is the same as in Figure 6.4. Upper left: Unperturbed Cartesian grid (GT1). Upper right: Prismatically extended simplex grid (GT2). Lower left: Cartesian grid with $\delta^{3/2}$ perturbations (GT3). Lower right: Cartesian grid with $\delta$ perturbations (GT4).

## 6.3 Poroelastic materials

Turning now to the full poroelastic system, we fix $\lambda = \vartheta = 1$, and vary the permeability according to $\kappa = \{1, 10^{-2}, 10^{-4}, 0\}$, exploring the limit of a vanishing permeability coefficient in Darcy's law. The Biot modulus and fluid compressibility are set to $\vartheta = 1$ and $\eta_w = 0.01$. In this case we consider the manufactured solution

$$u = \sum_{i=1}^{3} e_i x_{i+1}(1 - x_{i+1}) x_{i-1}(1 - x_{i-1}) \sin(2\pi x_{i+1}),$$

$$r = \sum_{i=1}^{3} e_i x_i (1 - x_i) \sin(\pi x_{i+1}) \sin(\pi x_{i-1}) \ ,$$

$$p = \prod_{i=1}^{3} x_i(1 - x_i), \qquad w = \sin(2\pi x_1)\, x_2(1 - x_2)x_3(1 - x_3). \tag{6.4}$$

Again, homogeneous Dirichlet boundary conditions are assigned.

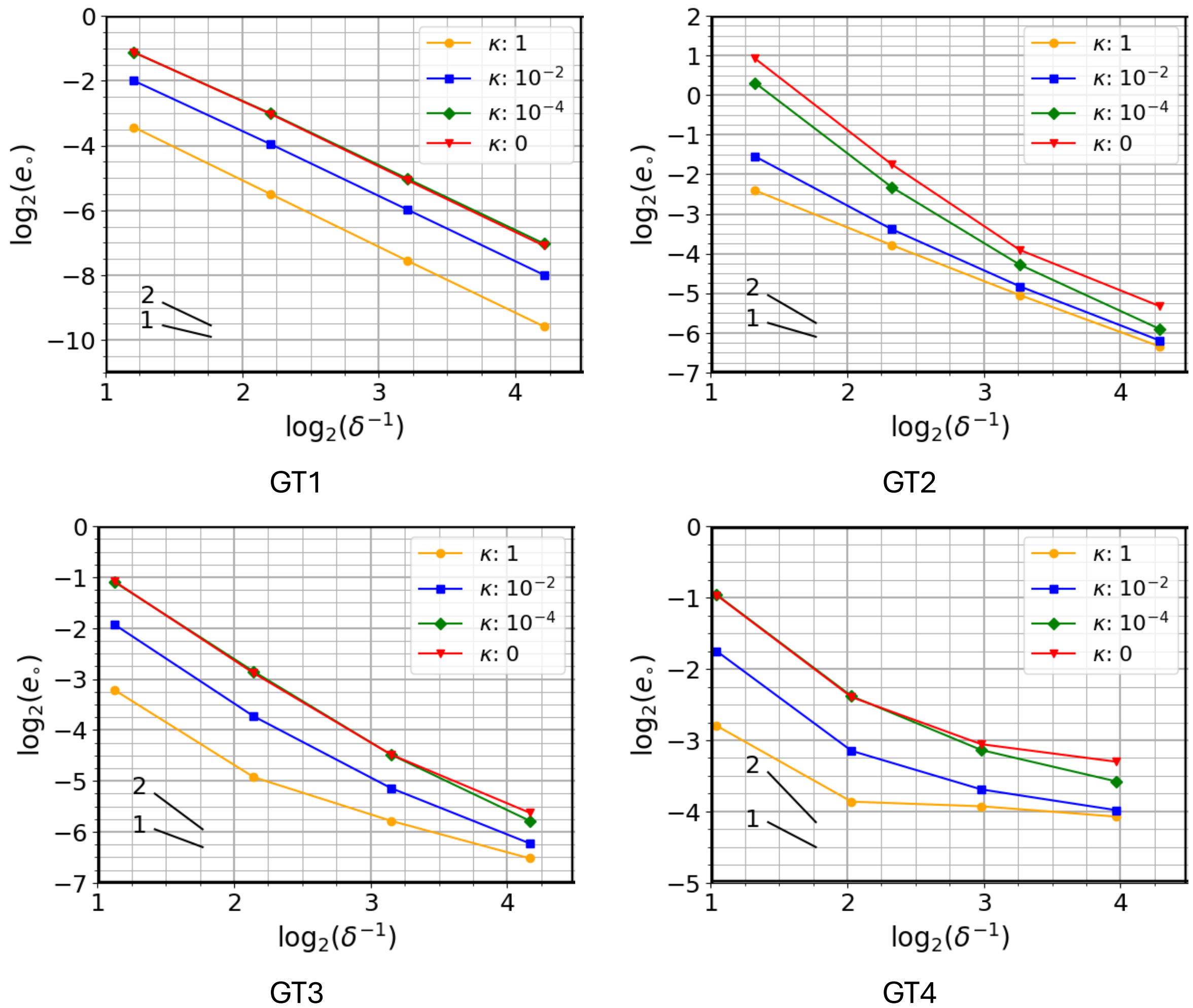

**Figure 6.6**: Convergence of TPSA vanishing parameter $\kappa$. Upper left: Unperturbed Cartesian grid (GT1). Upper left: Unperturbed Cartesian grid (GT1). Upper right: Prismatically extended simplex grid (GT2). Lower left: Cartesian grid with $\delta^{3/2}$ perturbations (GT3). Lower right: Cartesian grid with $\delta$ perturbations (GT4).

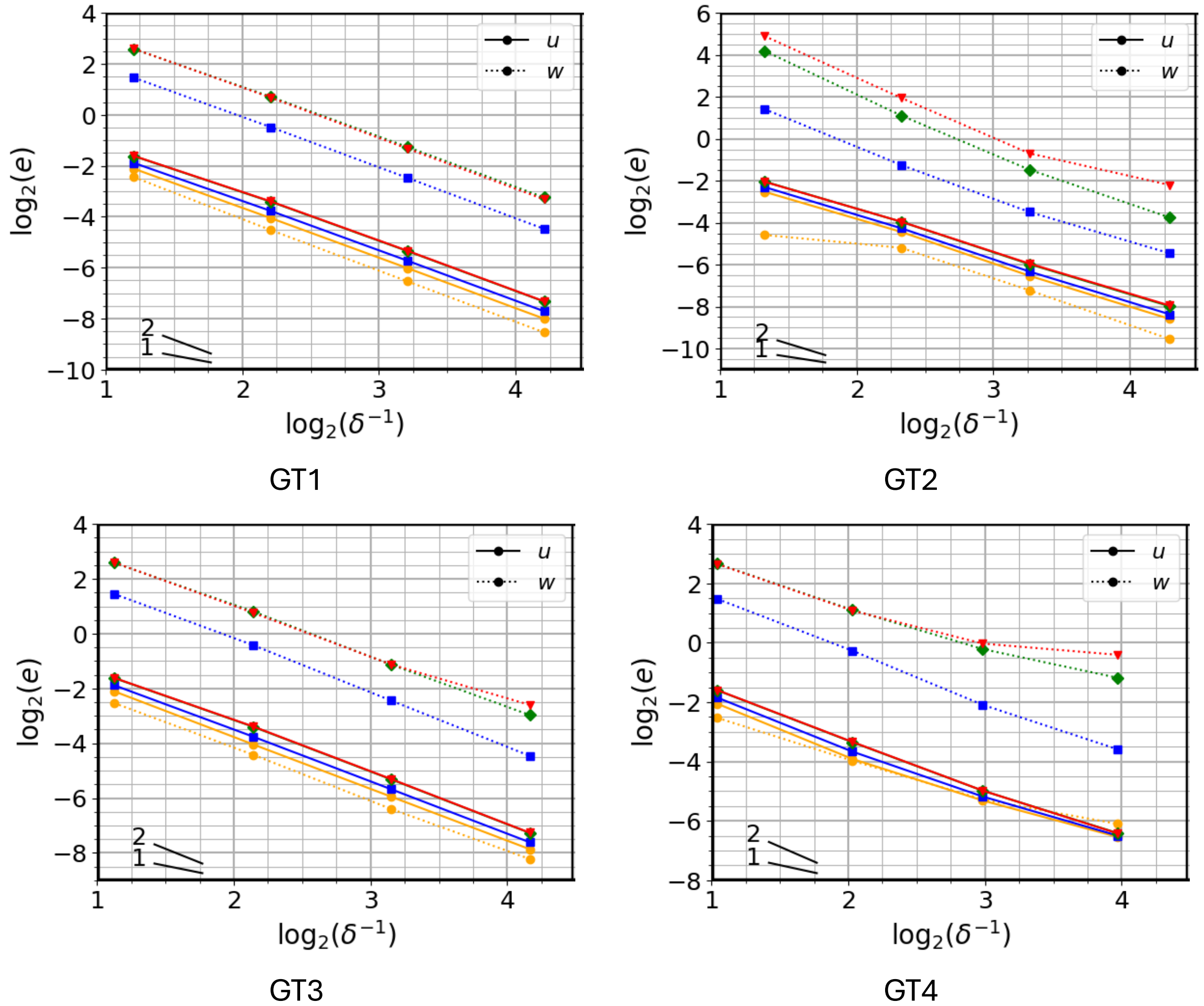

**Figure 6.7:** Convergence of the displacement and fluid pressure for a poroelastic material for a vanishing parameter $\kappa$. The color coding is the same as in Figure 6.6. Upper left: Unperturbed Cartesian grid (GT1). Upper right: Prismatically extended simplex grid (GT2). Lower left: Cartesian grid with $\delta^{3/2}$ perturbations (GT3). Lower right: Cartesian grid with $\delta$ perturbations (GT4).

The relative errors for our four classes of grid types are shown in Figure 6.6. As expected, the plots again confirm the stability of the method for all grid types and values of $\kappa$. Furthermore, as in the previous cases, the solution is second order convergent on GT1, and also converges on GT2 and GT3. On the most irregular grid, GT4, the error remains independent of $\delta$. The convergence of the primary variables $u$ and $w$ is shown in Figure 6.7. The displacement displays second order convergence on GT1, GT2 and GT3. The fluid pressure $w$ converges at second order when $\kappa > 0$, but falls to first order for $\kappa = 0$. This is expected, as the regularity of the solution drops from $H^1$ to $L^2$ in the absence of permeability, as reflected in the definition of the norm, equation (5.2). Overall, the primary variables again have higher convergence rates on GT2 and GT3 than has the total error $e_\circ$.

## 7. Concluding remarks

We have proposed and analyzed a new finite volume spatial discretization for the linearized elasticity, Stokes, and poromechanical equations, named TPSA. The discretization has the advantage that it does not use any form of dual or staggered grid, and that all variables are co-located using essentially the same spatial operators having the minimal sparsity pattern possible. The discretization is a generalization of the popular Two-Point Flux Approximation (TPFA) for flow in porous media, and as such, is expected to have similar strengths and weaknesses.

Both the theoretical and numerical analysis supports that the TPSA inherits the general properties of TPFA. In particular, the method has very strong stability properties, being robust on very general grids and for all relevant degeneracies of material parameters. On general grids, this comes at the cost of the consistency of the method, and indeed, convergence can only be expected for reasonably nice grids, with a so-called "face-orthogonal" property.

While this may seem somewhat restrictive, we emphasize the popularity of the TPFA method, which is the *de facto* standard discretization in all industrial computations for multi-face flow in porous media [43]. This suggests that the industry is well equipped to design grids for which the limitations of the TPFA and TPSA discretizations are manageable.

## Acknowledgements

The authors thank Wietse M. Boon and Omar Duran for fruitful discussions on the topic of Cosserat materials over several years. The work of EK was financed in part by Norwegian Research Council Grant 308733. The work of JMN was financed in part through his affiliation with NORCE and their Norwegian Research Council Grant 336294 (with support from Equinor, Norske Shell and Wintershall DEA Norge).

# Appendix A1: Calculus operators

Herein, we give an overview of the notational conventions used in the manuscript.

Throughout the paper, we consider a simply connected polygonal/polyhedral domain $\Omega \in \mathbb{R}^N$, where $N = 3$. The reduction to $N = 2$ is discussed in Appendix A4.

Throughout the paper, we adhere as much as possible to “matrix-vector” calculus, limiting the explicit use of indexing. This is as opposed to a more general tensor-based or exterior calculus-based exposition. We justify this as a compromise between simplicity of exposition and generality, which should be appropriate in the context of most numerical implementations.

We will denote primary variables such as displacement, solid rotation, solid pressure and fluid pressure by lower-case Latin letters, independent of whether they are scalars or vectors. We will denote secondary variables such as stresses and flux by lower case Greek letters, independent of whether they are vectors or matrices. We will also employ Greek letters for the (scalar) parameters of the problem. These are always considered functions of space. We reserve the use of capital letters for operators. When discussing discrete equations, we will use boldface to indicate a vector of discrete variables (such as e.g. the vector of cell-pressures).

Finally, we will use the relation $a \lesssim b$ (and $\gtrsim$) to imply that $a \leq Cb$, where the constant $C$ is robust in the sense of Definition 1.6.

## A1.1 Vector calculus conventions and definitions

Throughout this sub section, recall that $N = 3$, and let $p, q \in \mathbb{R}$ be scalars, $u, v, r \in \mathbb{R}^N$ be vectors with components $u_i$, $v_i$ and $r_i$, and $\sigma, \chi \in \mathbb{R}^{N\times N}$ be matrices with components $\sigma_{i,j}$ and $\chi_{i,j}$.

The matrix-vector product is defined as the vector $v = \sigma u$ such that $v_i = \sum_j \sigma_{i,j} u_j$. Similarly, the vector and matrix inner product is defined as $u \cdot v = \sum_i u_i v_i$ and $\sigma : \chi = \sum_{i,j} \sigma_{i,j} \chi_{i,j}$ .

When variables are functions over $\Omega$, we denote inner products as:

$$(p, q) \equiv \int_\Omega pq\, dV, \qquad (u, v) \equiv \int_\Omega u \cdot v\, dV, \qquad (\sigma, \chi) \equiv \int_\Omega \sigma \cdot \chi\, dV, \tag{A1.1}$$

For sufficiently smooth vector function $u$ and a matrix function $\sigma$, the divergence is defined as:

$$(\nabla \cdot u)(x) \equiv \lim_{\epsilon\to 0} \left|B_{x,\epsilon}\right|^{-1} \int_{\partial B_{x,\epsilon}} u \cdot n\, dA \qquad \text{and} \qquad (\nabla \cdot \sigma)(x) \equiv \lim_{\epsilon\to 0} \left|B_{x,\epsilon}\right|^{-1} \int_{\partial B_{x,\epsilon}} \sigma n\, dA. \tag{A1.2}$$

Here $B_{x,\epsilon}$ is an $N$-dimensional ball centered on $x$ of radius $\epsilon$. We define the gradient as the negative adjoint of the divergence, thus for scalar and vector functions $p$ and $v$,

$$(\nabla p, u) = -(p, \nabla \cdot u) \qquad \text{and} \qquad (\nabla v, \sigma) = -(v, \nabla \cdot \sigma), \tag{A1.3}$$

required to hold for all $C_0^\infty$ vector and matrix functions $u$ and $\sigma$.

The identity matrix $I$ is considered as the unique matrix such that $Iu = u$ for all $u$. Of importance in this work is the operator $S : \mathbb{R}^{N\times N} \to \mathbb{R}^N$, which intuitively measures the asymmetry of a matrix:

$$(S\sigma)_i \equiv \sigma_{i-1,i+1} - \sigma_{i+1,i-1}. \tag{A1.4}$$

Here and in the following indexes are always understood modulo $N$. The adjoint of $S$, denoted $S^* : \mathbb{R}^N \to \mathbb{R}^{N\times N}$, which satisfies, for all $\sigma$:

$$\sigma : S^* r = S\sigma \cdot r. \tag{A1.5}$$

The skew-symmetric matrix obtained from $S^*r$ has the explicit expression:

$$S^*r = \begin{pmatrix} 0 & -r_3 & r_2 \\ r_3 & 0 & -r_1 \\ -r_2 & r_1 & 0 \end{pmatrix}. \tag{A1.6}$$

The fact that $S^*$ is the right inverse of $S$ (up to a factor 2) will be used frequently, thus we summarize:

$$SS^*r = 2r \qquad \text{while} \qquad S^*S\sigma = \sigma - \sigma^T. \tag{A1.7}$$

The operator $S^*$ is itself skew-self adjoint:

$$(S^*r)u = -(S^*u)r. \tag{A1.8}$$

We will frequently need the expression in the above parenthesis, in particular as applied to a normal vector $n$. Recognizing that $S^*n$ is a rotation matrix around the axis provided by $n$, leads to the shorthand $R^n = S^*n$. Finally, we note the useful relationship

$$S\nabla u = -\nabla \cdot S^*u. \tag{A1.9}$$

## A1.2 Continuous and discrete function spaces

For the continuous variables, we employ standard Hilbert spaces of functional analysis. Based on the standard inner products (see equation (A1.1)), we denote the $\gamma$-weighted $L^2$ norm for scalars, vectors and matrices as:

$$\|u\|_\gamma = \sqrt{(\gamma u, u)}. \tag{A1.10a}$$

As a convention, we will omit the subscript when no weight is considered, i.e. $\|u\| = \|u\|_{\gamma=1}$. Additionally, we will need the seminorm:

$$|u|_\gamma = \|u - \bar{u}\|_\gamma. \tag{A1.10b}$$

where $\bar{u} = |\Omega|^{-1}(1, u)$ is the mean value of $u$ over $\Omega$.

We will employ discrete inner products for both cell-variables $\boldsymbol{u}, \boldsymbol{v} \in V^{|I_\mathcal{T}|}$ and face-variables $\boldsymbol{\sigma}, \boldsymbol{\psi} \in V^{|I_\mathcal{F}|}$, where $V = \mathbb{R}^p$ the (local) vector space of the variable for some $p \geq 0$. Thus

$$(\boldsymbol{u}, \boldsymbol{v}) = \textstyle\sum_{i\in I_\mathcal{T}} u_i \cdot v_i \qquad \text{and} \qquad (\boldsymbol{\sigma}, \boldsymbol{\psi}) = \textstyle\sum_{i\in I_\mathcal{T}} \sigma_i \cdot \psi_i. \tag{A1.11}$$

The corresponding weighted norms are defined respectively as:

$$\|\boldsymbol{u}\|_\gamma = \sqrt{\textstyle\sum_{i\in I_{\mathcal{T}_\delta}} |\omega_i| \gamma_i u_i \cdot u_i} \qquad \text{and} \qquad \|\boldsymbol{\psi}\|_\gamma = \sqrt{\textstyle\sum_{k\in I_{\mathcal{F}_\delta}} \frac{|\varsigma_k|\delta_k}{N} \gamma_k \psi_k \cdot \psi_k}. \tag{A1.12}$$

As in the continuous case, we will omit subscripts if the weight is unity. Note that for discrete variables, we weight the norms by volumetric quantities, but not the inner products.

To relate the discrete and continuous norms, we define the $\pi_\omega$ as the $L^2$ projection on to the piecewise constants $P_0(\mathcal{T})$ and $\pi_\varsigma$ as the $L^2$ projection of the normal component onto the piecewise constants $P_0(\mathcal{F})$, e.g.:

$$(\pi_\omega u)_i = |\omega_i|^{-1} \int_{\omega_i} u \, dV \quad \text{and} \quad (\pi_\varsigma \psi)_k = |\varsigma_k|^{-1} \int_{\omega_k} \psi \cdot n \, dV. \tag{A1.13}$$

Then it holds that [3]:

$$\|u\|_\gamma \leq \|\pi_\omega u\|_\gamma \quad \text{and} \quad \|\psi\|_\gamma \lesssim \|\pi_\varsigma \psi\|_\gamma. \tag{A1.14}$$

We also note that whenever there is some part of the boundary with zero Dirichlet data, both continuous and discrete Poincaré inequalities hold [23]:

$$\|u\| \lesssim \|\nabla u\| \quad \text{and} \quad \|\boldsymbol{u}\| \lesssim \|\boldsymbol{\delta}^{-1} \Delta^* \boldsymbol{u}\|. \tag{A1.15}$$

A stronger result also holds in the continuous case, known as Korn's inequality, which states that subject to an appropriate side constraint (such as e.g. a non-trivial section of the boundary with $u = 0$), there exists a $0 < C_K$ such that [44]:

$$C_K \|\nabla u\|_\mu \leq \left\| \frac{\nabla u + \nabla u^T}{2} \right\|_\mu. \tag{A1.16}$$

# Appendix A2: Derivation of TPSA coefficients

In this section, we provide the detailed derivation of the expressions appearing in Section 3.1.

We approximate the normal stresses of an internal face $\varsigma_k$ based on the variables defined in the two cells $\omega_i$ and $\omega_j$ for $\{i, j\} \in \mathcal{N}_k^*$. Throughout the appendix, we assume without loss of generality that the indexes are ordered such that $n_{k,i} = n_k = -n_{k,j}$. We first recall that the definition of the numerical flux implies that:

$$\psi_k = \int_{\varsigma_k} \psi \cdot n_k \, dA \approx |\varsigma_k| \psi \cdot n_k. \tag{A2.1}$$

The key point is therefore to approximate $\psi \cdot n_k$, based on the constitutive law given in equation (3.2). Seen from cell $\omega_i$, we therefore obtain the approximations:

$$|\varsigma_k|^{-1} \sigma_k = 2\mu_i (\nabla u \cdot n_k)_{k,i} + S^* r_{k,i} n_k + p_{k,i} n_k, \tag{A2.2a}$$

$$|\varsigma_k|^{-1} \tau_k = S^* u_{k,i} \cdot n_k, \tag{A2.2b}$$

$$|\varsigma_k|^{-1}v_k = u_{k,i} \cdot n_k. \tag{A2.2c}$$

The normal derivatives are naturally approximated based on the difference between the cell-center value and some value $z_k$ of the primary variable at the face $\varsigma_k$. Concretely, we use the standard difference formula:

$$(\nabla z \cdot n_k)_{k,i} \approx \frac{z_k - z_i}{\delta_k^i} = \delta_k^{-i}(z_k - z_i). \tag{A2.3}$$

Here we introduce the useful shorthand notation $\delta_k^{-i} = \left(\delta_k^i\right)^{-1}$.

In equations (A2.2), we also see three algebraic expressions $(2\mu S^* r)_{k,i}$, $p_{k,i}$ and $(2\mu S^* u)_{k,i}$, that need to be determined, in the sense of the expression evaluated on face $k$ as seen from cell $i$. From the analysis of the continuous problem in Section 5.1, we always expect displacement to have $H^1$ regularity, while we only expect $L^2$ regularity of solid pressure. This motivates setting

$$u_{k,i} = u_k, \qquad \text{and} \qquad p_{k,i} = p_i. \tag{A2.4a}$$

For the rotations, we also formally only expect $L^2$ regularity. On the other hand, from a modeling perspective, higher regularity of rotations may often be natural [45, 19]. We therefore allow for this possibility, and retain this flexibility, denoted by $\alpha$ and $\tilde{\alpha} = 1 - \alpha$ in the derivation:

$$r_{k,i} = \alpha r_i + \tilde{\alpha} r_k. \tag{A2.4b}$$

In summary, equations (2.2) are represented as:

$$|\varsigma_k|^{-1}\sigma_k = 2\mu_i \delta_k^{-i}(u_k - u_i) + S^*(\alpha r_i + \tilde{\alpha} r_k) n_k + p_i n_k, \tag{A2.5a}$$

$$|\varsigma_k|^{-1}\tau_k = S^* u_k \cdot n_k, \tag{A2.5b}$$

$$|\varsigma_k|^{-1}v_k = u_k \cdot n_k. \tag{A2.5c}$$

## A2.1 Discretization stencils for internal faces

For internal faces, Equations (A2.5a) can be restated as seen from $\omega_j$:

$$|\varsigma_k|^{-1}\sigma_k = 2\mu_j \delta_k^{-i}\left(u_j - u_k\right) + S^*\left(\alpha r_j + \tilde{\alpha} r_k\right) n_k + p_j n_k. \tag{A2.6a}$$

Equation (A2.5c) does not depend on $i$, and is thus identical from both sides. Equations (A2.5-A2.6) provide sufficient constraints to eliminate the intermediate primary variables $u_k$ and $r_k$, and obtain an expression for the numerical flux only in terms of the cell-variables $z_i$ and $z_j$. To show this, we use the reorder the product terms to get the variables to the right, using the relations (see Appendix A1 and equation (3.6)):

$$u \cdot n = nu, \qquad (S^* u) \cdot n = -R^n u, \qquad (S^* r) n = -R^n r. \tag{A2.7}$$

Using the notation $R_k^n \equiv R^{n_k} = S^* n_k$, we collect equations (A2.5-A2.6) in a linear system as follows:

$$\begin{pmatrix} 1 & & & -2\mu_i\delta_k^{-i} \\ & 1 & & R_k^n \\ & & 1 & -n_k \\ 1 & & & 2\mu_j\delta_k^{-j} \end{pmatrix}\begin{pmatrix} \sigma_k \\ \tau_k \\ \upsilon_k \\ |\varsigma_k|u_k \end{pmatrix} = |\varsigma_k| \begin{pmatrix} -2\mu_i\delta_k^{-i} & -\alpha R_k^n & n_k & & & & -\tilde{\alpha}R_k^n \\ & & & & & & \\ & & & & & & \\ & & & 2\mu_j\delta_k^{-j} & -\alpha R_k^n & n_k & -\tilde{\alpha}R_k^n \end{pmatrix}\begin{pmatrix} u_i \\ r_i \\ p_i \\ u_j \\ r_j \\ p_j \\ r_k \end{pmatrix}. \tag{A2.8}$$

The above equations have the following structure, (in the following equations through equation (A2.18) below, the Latin letters $a$ through $h$ have no relation to their use in the main part of the manuscript):

$$\begin{pmatrix} 1 & & & -a \\ & 1 & & b \\ & & 1 & -d \\ 1 & & & e \end{pmatrix}\begin{pmatrix} \sigma_k \\ \tau_k \\ \upsilon_k \\ |\varsigma_k|u_k \end{pmatrix} = |\varsigma_k| \begin{pmatrix} -a & -\alpha b & d & & & & -\tilde{\alpha}b \\ & & & & & & \\ & & & & & & \\ & & & e & -\alpha b & d & -\tilde{\alpha}b \end{pmatrix}\begin{pmatrix} u_i \\ r_i \\ p_i \\ u_j \\ r_j \\ p_j \\ r_k \end{pmatrix}. \tag{A2.9}$$

We solve this system directly, to obtain:

$$\begin{pmatrix} \sigma_k \\ \tau_k \\ \upsilon_k \\ |\varsigma_k|u_k \end{pmatrix} = |\varsigma_k| \begin{pmatrix} -a + ac_1a & -\alpha b + ac_1\alpha b & d - ac_1d & ac_1e & -ac_1\alpha b & ac_1d & -\tilde{\alpha}b \\ -bc_1a & -bc_1\alpha b & bc_1d & -bc_1e & bc_1\alpha b & -bc_1d & 0 \\ dc_1a & dc_1\alpha b & -dc_1d & dc_1e & -dc_1\alpha b & dc_1d & 0 \\ c_1a & c_1\alpha b & -c_1d & c_1e & -c_1\alpha b & c_1d & 0 \end{pmatrix}\begin{pmatrix} u_i \\ r_i \\ p_i \\ u_j \\ r_j \\ p_j \\ r_k \end{pmatrix}, \tag{A2.10}$$

where $c_1 = (a+e)^{-1}$.

To simplify this expression, the following relationships are useful:
*Algebraic identity:*

$$(1 - ac_1) = c_1e. \tag{A2.11}$$

*Geometric identities*:

$$dc_1d = c_1n_kn_k = c_1, \qquad bc_1d = c_1R_k^n n_k = 0 = dc_1f\,. \tag{A2.12}$$

Thus from (A2.10) we obtain:

$$\begin{pmatrix} \sigma_k \\ \tau_k \\ \upsilon_k \\ |\varsigma_k| u_k \end{pmatrix} = |\varsigma_k| \begin{pmatrix} -ec_1a & -ec_1\alpha b & ec_1d & ac_1e & ac_1\alpha b & ac_1d & \tilde{\alpha}b \\ -bc_1a & -bc_1\alpha b & 0 & -bc_1e & bc_1\alpha b & 0 & 0 \\ dc_1a & 0 & -c_1 & dc_1e & 0 & c_1 & 0 \\ c_1a & c_1\alpha b & -c_1d & c_1e & c_1\alpha b & c_1d & 0 \end{pmatrix} \begin{pmatrix} u_i \\ r_i \\ p_i \\ u_j \\ r_j \\ p_j \\ r_k \end{pmatrix}. \tag{A2.13}$$

We recognize in equation (A2.13) that the first three columns are very similar to the next three columns. Indeed, by comparison to Section 2.1, we note that that e.g.

$$u_i - u_j = \Delta_k^* \boldsymbol{u}, \tag{A2.14a}$$

$$c_1\left(au_i + eu_j\right) = \frac{\mu_i \delta_k^{-i} u_i + \mu_j \delta_k^{-j} u_j}{\mu_i \delta_k^{-i} + \mu_j \delta_k^{-j}} = \Xi_k \boldsymbol{u}, \tag{A2.14b}$$

$$bc_1\left(er_{ri} + ar_j\right) = b\,\frac{\mu_j \delta_k^{-j} r_i + \mu_i \delta_k^{-i} r_j}{\mu_i \delta_k^{-i} + \mu_j \delta_k^{-j}} = b\tilde{\Xi}_k \boldsymbol{r}, \tag{A2.14c}$$

$$dc_1\left(ep_i + ap_j\right) = d\,\frac{\mu_j \delta_k^{-j} p_i + \mu_i \delta_k^{-i} u_j}{\mu_i \delta_k^{-i} + \mu_j \delta_k^{-j}} = d\tilde{\Xi}_k \boldsymbol{p}. \tag{A2.14d}$$

Moreover, we realize that the face rotation $r_k$ is unspecified. Motivated by (A2.14c), we avoid introducing additional averaging operators and define

$$r_k = \tilde{\Xi}_k \boldsymbol{r}. \tag{A2.15}$$

With these substitutions, equation (A2.13) simplifies further as:

$$\begin{pmatrix} \sigma_k \\ \tau_k \\ \upsilon_k \\ |\varsigma_k| u_k \end{pmatrix} = |\varsigma_k| \begin{pmatrix} -ec_1a & & & & -b & d \\ & -bc_1\alpha b & & -b & & \\ & & -c_1 & d & & \\ & c_1\alpha b & -c_1d & 1 & & \end{pmatrix} \begin{pmatrix} \Delta_k^* \boldsymbol{u} \\ \Delta_k^* \boldsymbol{r} \\ \Delta_k^* \boldsymbol{p} \\ \Xi_k \boldsymbol{u} \\ \tilde{\Xi}_k \boldsymbol{r} \\ \tilde{\Xi}_k \boldsymbol{p} \end{pmatrix}. \tag{A2.16}$$

The first three lines of Equation (A2.16) represent the desired numerical fluxes. Before restating these by substituting the definitions of the compound variables, we recall the mean values defined in equation (3.5):

*$\mu$-weighted distance:*

$$c_1 = \left(2\mu_i \delta_k^{-i} + 2\mu_j \delta_k^{-j}\right)^{-1} = \delta_k^{\mu}. \tag{A2.17}$$

*Harmonic means:*

$$ac_1e = 2\,\frac{\mu_i \delta_k^{-i} \mu_j \delta_k^{-j}}{\mu_i \delta_k^{-i} + \mu_j \delta_k^{-j}} = \delta_k^{-1} 2\bar{\mu}_k. \tag{A2.18}$$

From these identifications, we obtain the numerical TPSA fluxes:

$$\begin{pmatrix}\sigma_k\\ \tau_k\\ \upsilon_k\end{pmatrix} = |\varsigma_k| \begin{pmatrix} -\delta_k^{-1}2\bar{\mu}_k & & & & -R_k^n & n_k \\ & -\alpha\delta_k^{\mu}(R_k^n)^2 & & -R_k^n & & \\ & & -\delta_k^{\mu} & n_k & & \end{pmatrix} \begin{pmatrix}\Delta_k^*\boldsymbol{u}\\ \Delta_k^*\boldsymbol{r}\\ \Delta_k^*\boldsymbol{p}\\ \Xi_k\boldsymbol{u}\\ \tilde{\Xi}_k\boldsymbol{r}\\ \tilde{\Xi}_k\boldsymbol{p}\end{pmatrix}. \quad \text{(A2.19)}$$

The operator $(R_k^n)^2$ is the negative of the projection onto the plane orthogonal to $n_k$. As such, the term $\alpha\delta_k^{\mu}(R_k^n)^2$ has the character of a negative Laplacian, and is destabilizing in the method. This motivates setting $\alpha \leq 0$, and indeed, we adopt the simplest choice $\alpha = 0$ for internal faces. Thus, we obtain in terms of the primary cell-centered variables $\boldsymbol{u}$, $\boldsymbol{r}$ and $\boldsymbol{p}$, which is reported in equation (3.7) in the main part of the manuscript.

## A2.2 Boundary faces

Recall that for boundary faces $\varsigma_k \in \mathcal{F} \cap \partial\Omega$, it holds that $N_k^* = \{i, j\}$, where without loss of generality, we can assume that $i \in I_{\mathcal{T}}$ and $j = k \in I_{\mathcal{B}}$. The boundary conditions are stated in equations (1.16), and when specified onto a face become a Robin BC for displacement with length scale $\delta_k^{u,k} = b_k^u$:

$$\delta_k^{u,k}\left(\sigma_k - |\varsigma_k| g_{\sigma,k}^u\right) = 2\mu\left(|\varsigma_k| g_{u,k}^u - |\varsigma_k| u_k\right). \quad \text{(A2.21)}$$

In view of this, equations (A2.8) already provide the structure necessary to consider general boundary conditions. At the boundary, $r_k$ is not defined, and we only carry out the derivation for $\alpha = 1$:

$$\begin{pmatrix} 1 & & & -2\mu_i\delta_k^{-i} \\ & 1 & & R_k^n \\ & & 1 & -n_k \\ 1 & & & 2\mu_k\delta_k^{u,-k} \end{pmatrix} \begin{pmatrix}\sigma_k\\ \tau_k\\ \upsilon_k\\ |\varsigma_k| u_k\end{pmatrix} = |\varsigma_k| \begin{pmatrix} -2\mu_i\delta_k^{-i} & -R_k^n & n_k \\ & & \\ & & \\ & & \end{pmatrix} \begin{pmatrix}u_i\\ r_i\\ p_i\end{pmatrix} + |\varsigma_k| \begin{pmatrix} & \\ & \\ & \\ 1 & 2\mu_k\delta_k^{u,-k} \end{pmatrix} \begin{pmatrix} g_{\sigma,k}^u \\ g_{u,k}^u \end{pmatrix}. \quad \text{(A2.22)}$$

The same derivation as in Section A2.1 applies, and we obtain:

$$\begin{pmatrix}\sigma_k\\ \tau_k\\ \upsilon_k\\ |\varsigma_k| u_k\end{pmatrix} = |\varsigma_k| \begin{pmatrix} -ec_1a & -ec_1b & ec_1d \\ -bc_1a & -bc_1b & 0 \\ dc_1a & 0 & -c_1 \\ c_1a & c_1b & -c_1d \end{pmatrix} \begin{pmatrix}u_i\\ r_i\\ p_i\end{pmatrix} + |\varsigma_k| \begin{pmatrix} ac_1 & ac_1e \\ -bc_1 & -bc_1e \\ dc_1 & dc_1e \\ c_1 & c_1e \end{pmatrix} \begin{pmatrix} g_{\sigma,k}^u \\ g_{u,k}^u \end{pmatrix}. \quad \text{(A2.23)}$$

We note that

$$c_1 g_{\sigma,k}^u + c_1 e g_{u,k}^u = \delta_k^{\mu}\left(g_{\sigma,k}^u + 2\mu_k\delta_k^{u,-k} g_{u,k}^u\right) = -\delta_k^{\mu}\Delta_{\mathcal{B},k}^* g_{\sigma,k}^u + \Xi_{\mathcal{B},k} g_{u,k}^u. \quad \text{(A2.24)}$$

Thus, we obtain the general expression for the boundary conditions as:

$$\begin{pmatrix}\sigma_k\\ \tau_k\\ \upsilon_k\end{pmatrix} = |\varsigma_k|\begin{pmatrix}-\delta_k^{-1}2\bar{\mu}_k\Delta^*_{\mathcal{T},k} & -R_k^n\tilde{\Xi}_{\mathcal{T},k} & n_k\tilde{\Xi}_{\mathcal{T},k}\\ -R_k^n\Xi_{\mathcal{T},k} & -\delta_k^{\mu}(R_k^n)^2\Delta^*_{\mathcal{T},k} & \\ n_k\Xi_{\mathcal{T},k} & & -\delta_k^{\mu}\Delta^*_{\mathcal{T},k}\end{pmatrix}\begin{pmatrix}\boldsymbol{u}\\ \boldsymbol{r}\\ \boldsymbol{p}\end{pmatrix} + |\varsigma_k|\begin{pmatrix}\tilde{\Xi}_{\mathcal{B},k} & -\delta_k^{-1}2\bar{\mu}_k\Delta^*_{\mathcal{B},k}\\ \delta_k^{\mu}R_k^n\Delta^*_{\mathcal{B},k} & -R_k^n\Xi_{\mathcal{B},k}\\ -\delta_k^{\mu}n_k\Delta^*_{\mathcal{B},k} & n_k\Xi_{\mathcal{B},k}\end{pmatrix}\begin{pmatrix}\boldsymbol{g}_\sigma^u\\ \boldsymbol{g}_u^u\end{pmatrix}. \qquad \text{(A2.25)}$$

In this equation, we emphasize that both $\delta_k$ and $\Xi_k$ depend on $b^u$. It is often useful to identify the displacement of the boundary (e.g. for contact mechanics), and from equation (A2.25) we also obtain:

$$u_k = \Xi_{\mathcal{T},k}\boldsymbol{u} + \delta_k^{\mu}R_k^n\Delta^*_{\mathcal{T},k}\boldsymbol{r} - \delta_k^{\mu}n_k\Delta^*_{\mathcal{T},k}\boldsymbol{p} + \delta_k^{\mu}g^u_{\sigma,k} + \Xi_{\mathcal{B},k}g^u_{u,k}.$$

From Equation (A2.25), we can simplify the boundary conditions for the three most common cases:

- Dirichlet boundary conditions $b^u = 0$: In this case, $\left(\Xi_{\mathcal{B},k}\right)_k = 1$, while $\left(\Xi_{\mathcal{T},k}\right)_i = \delta_k^{\mu} = 0$, thus Equation (A2.28) simplifies to:
$$\begin{pmatrix}\sigma_k\\ \tau_k\\ \upsilon_k\end{pmatrix} = |\varsigma_k|\begin{pmatrix}-\delta_k^{-1}2\bar{\mu}_k\Delta^*_{\mathcal{T},k} & -R_k^n\tilde{\Xi}_{\mathcal{T},k} & n_k\tilde{\Xi}_{\mathcal{T},k}\\ & & \\ & & \end{pmatrix}\begin{pmatrix}\boldsymbol{u}\\ \boldsymbol{r}\\ \boldsymbol{p}\end{pmatrix} + |\varsigma_k|\begin{pmatrix}-\delta_k^{-1}2\bar{\mu}_k\Delta^*_{\mathcal{B},k}\\ -R_k^n\Xi_{\mathcal{B},k}\\ n_k\Xi_{\mathcal{B},k}\end{pmatrix}\boldsymbol{g}_u^u. \qquad \text{(A2.26)}$$
- Neumann boundary conditions: $(b^u)^{-1} = 0$: In this case $\left(\Xi_{\mathcal{B},k}\right)_k = \delta_k^{-1} = 0$, while $\delta_k^{\mu} = (2\mu_i)^{-1}\delta_k^i$, thus and Equation (A2.25) simplifies to:
$$\begin{pmatrix}\sigma_k\\ \tau_k\\ \upsilon_k\end{pmatrix} = |\varsigma_k|\begin{pmatrix} & & \\ -R_k^n\Xi_{\mathcal{T},k} & -\delta_k^{\mu}(R_k^n)^2\Delta^*_{\mathcal{T},k} & \\ n_k\Xi_{\mathcal{T},k} & & -\delta_k^{\mu}\Delta^*_{\mathcal{T},k}\end{pmatrix}\begin{pmatrix}\boldsymbol{u}\\ \boldsymbol{r}\\ \boldsymbol{p}\end{pmatrix} + |\varsigma_k|\begin{pmatrix}1\\ \delta_k^{\mu}R_k^n\Delta^*_{\mathcal{B},k}\\ -\delta_k^{\mu}n_k\Delta^*_{\mathcal{B},k}\end{pmatrix}\boldsymbol{g}_\sigma^u. \qquad \text{(A2.27)}$$
- For homogeneous s boundary conditions, $\boldsymbol{g}_u^u = \boldsymbol{g}_\sigma^u = 0$, general Robin boundary conditions can all be written as:
$$\begin{pmatrix}\sigma_k\\ \tau_k\\ \upsilon_k\end{pmatrix} = |\varsigma_k|\begin{pmatrix}-\delta_k^{-1}2\bar{\mu}_k\Delta^*_{\mathcal{T},k} & -R_k^n\tilde{\Xi}_{\mathcal{T},k} & n_k\tilde{\Xi}_{\mathcal{T},k}\\ -R_k^n\Xi_{\mathcal{T},k} & -\delta_k^{\mu}(R_k^n)^2\Delta^*_{\mathcal{T},k} & \\ n_k\Xi_{\mathcal{T},k} & & -\delta_k^{\mu}\Delta^*_{\mathcal{T},k}\end{pmatrix}\begin{pmatrix}\boldsymbol{u}\\ \boldsymbol{r}\\ \boldsymbol{p}\end{pmatrix}. \qquad \text{(A2.28)}$$
where we understand that the dependence on $b^u$ is fully captured by the definitions of the operators involved.

# Appendix A3: Proofs for Lemma 5.9

*Point 1.*: We calculate

$$\|\tilde{\Xi}\boldsymbol{p}\|^2_{\bar{\mu}^{-1}} = \sum_{k\in I_{\mathcal{F}}}\bar{\mu}_k^{-1}|\varsigma_k|\delta_k\left(\tilde{\Xi}_k\boldsymbol{p}\right)^2 \le \sum_{k\in I_{\mathcal{F}}}\bar{\mu}_k^{-1}|\varsigma_k|\delta_k\sum_{i\in N_k^*}\tilde{\Xi}_{k,i}p_i^2. \qquad \text{(A3.1)}$$

Now since (recalling the shorthand $\delta_k^{-i} = \left(\delta_k^i\right)^{-1}$):

$$\delta_k\tilde{\Xi}_{k,i} = \delta_k\frac{\mu_j\delta_k^{-j}}{\mu_i\delta_k^{-i}+\mu_j\delta_k^{-j}} = \frac{\bar{\mu}_k}{\mu_i\delta_k^{-i}}. \qquad \text{(A3.2)}$$

We obtain:

$$\left\|\tilde{\Xi}\boldsymbol{p}\right\|_{\overline{\mu}^{-1}}^{2} \leq \sum_{k\in I_{\mathcal{F}}} \frac{|\varsigma_k|}{N} \sum_{i\in\mathcal{N}_k^*} \delta_k^i \mu_i^{-1} p_i^2 = \frac{1}{N} \sum_{i\in I_{\mathcal{T}}} \mu_i^{-1} p_i^2 \sum_{k\in\mathcal{N}_i} |\varsigma_k| \delta_k^i. \tag{A3.3}$$

Now since it follows from the divergence theorem that:

$$\sum_{k\in\mathcal{N}_i} |\varsigma_k| \delta_k^i = \int_{\partial\omega_i} (x - x_i) \cdot n \, dA = \int_{\omega_i} \nabla \cdot x \, dV = N|\omega_i|. \tag{A3.4}$$

We can conclude the claim:

$$\left\|\tilde{\Xi}\boldsymbol{p}\right\|_{\overline{\mu}^{-1}}^{2} \leq \sum_{i\in I_{\mathcal{T}}} \mu_i^{-1} p_i^2 = \|\boldsymbol{p}\|_{\mu^{-1}}^2. \tag{A3.5}$$

*Point 2.*: We note that $-R^n R^n$ is a projection onto the plane orthogonal to $n$, and thus $(-R^n R^n)^2 = -R^n R^n$. A direct calculation then gives:

$$(\Delta\boldsymbol{u})_k = -R_k^n R_k^n (\Delta\boldsymbol{u})_k + n_k (n_k \cdot (\Delta\boldsymbol{u})_k). \tag{A3.6}$$

From which the identity follows.

*Point 3*: By the definition of adjoints:

$$(\Delta|\boldsymbol{\varsigma}|\boldsymbol{n}\Xi\boldsymbol{u}, \boldsymbol{p}) = (\Xi^*|\boldsymbol{\varsigma}|\boldsymbol{n}\Delta^*\boldsymbol{p}, \boldsymbol{u}), \tag{A3.7}$$

thus

$$(\Delta|\boldsymbol{\varsigma}|\boldsymbol{n}\Xi\boldsymbol{u}, \boldsymbol{p}) + \left(\Delta|\boldsymbol{\varsigma}|\boldsymbol{n}\tilde{\Xi}\boldsymbol{p}, \boldsymbol{u}\right) = \left(\left(\Delta|\boldsymbol{\varsigma}|\boldsymbol{n}\tilde{\Xi} + \Xi^*|\boldsymbol{\varsigma}|\boldsymbol{n}\Delta^*\right)\boldsymbol{p}, \boldsymbol{u}\right). \tag{A3.8}$$

However,

$$\begin{aligned}\left(\left(\Delta|\boldsymbol{\varsigma}|\boldsymbol{n}\tilde{\Xi} + \Xi^*|\boldsymbol{\varsigma}|\boldsymbol{n}\Delta^*\right)\boldsymbol{p}\right)_i &= \frac{1}{N} \sum_{k\in N_i} \left( |\varsigma_k| \Delta_{i,k} n_k \sum_{j\in N_k^*} \tilde{\Xi}_{k,j} p_j + \Xi_{k,i} n_k \sum_{j\in N_k^*} \Delta_{j,k} p_j \right) \\ &= \frac{1}{N} \sum_{k\in N_i} \left( |\varsigma_k| \Delta_{i,k} n_k \sum_{j\in N_k^*} \left( \tilde{\Xi}_{k,j} + \Xi_{k,i} \frac{\Delta_{j,k}}{\Delta_{i,k}} \right) p_j \right) = \frac{1}{N} \sum_{k\in N_i} \left( |\varsigma_k| \Delta_{i,k} n_k p_i \right) = 0,\end{aligned} \tag{A3.9}$$

since

$$\tilde{\Xi}_{k,j} + \Xi_{k,i} \frac{\Delta_{j,k}}{\Delta_{i,k}} = \begin{cases} 1 & \text{if} \quad i = j \\ 0 & \text{if} \quad i \neq j \end{cases} \quad \text{and} \quad \sum_{k\in N_i} \left( |\varsigma_k| \Delta_{i,k} n_k \right) = 0. \tag{A3.10}$$

*Point 4.*: By the definition of adjoints and the fact that $(R^{\boldsymbol{n}})^* = -R^{\boldsymbol{n}}$:

$$(\Delta|\boldsymbol{\varsigma}|R^{\boldsymbol{n}}\Xi\boldsymbol{u}, \boldsymbol{r}) = -(\Xi^*|\boldsymbol{\varsigma}|R^{\boldsymbol{n}}\Delta^*\boldsymbol{r}, \boldsymbol{u}), \tag{A3.11}$$

thus

$$(\Delta|\boldsymbol{\varsigma}|R^{\boldsymbol{n}}\Xi\boldsymbol{u},\boldsymbol{r}) - (\Delta|\boldsymbol{\varsigma}|R^{\boldsymbol{n}}\tilde{\Xi}\boldsymbol{r},\boldsymbol{u}) = \left((\Delta|\boldsymbol{\varsigma}|R^{\boldsymbol{n}}\tilde{\Xi} + \Xi^*|\boldsymbol{\varsigma}|R^{\boldsymbol{n}}\Delta^*)\boldsymbol{r},\boldsymbol{u}\right). \tag{A3.12}$$

Proceeding as in point 3 concludes the proof, since $R^{n_k} = S^* n_k$ and $S^*$ is a linear operator.

# Appendix A4: Reduction to 2D

In this Appendix, we consider a simply connected polygonal/polyhedral domain $\Omega \in \mathbb{R}^N$, where $N \in 2$. We do not treat the extension to poroelasticity, since this is equivalent to the 3D case.

The 2D reduction of mechanics is achieved by recognizing that a 2D domain can always be realized as a slice of a 3D domain where symmetry is imposed in the third dimension. We therefore denote the embedding of a 2D domain into a 3D domain by a dot above the variable, and consider the so-called "plane strain" representation, where the (vector) displacement $u$ is zero in the third dimension, while the (now scalar) rotation $r$ measures in-plane rotations of the first two dimensions. Solid pressure $p$ remains scalar, such that:

$$\dot{u} \equiv \begin{pmatrix} u_1 \\ u_2 \\ 0 \end{pmatrix}, \qquad \dot{r} \equiv \begin{pmatrix} 0 \\ 0 \\ r \end{pmatrix}, \qquad \text{and} \qquad \dot{p} = p. \tag{A4.1}$$

This motivates introducing the operators $I_{||}$ and $I_{\perp}$, encoding the parallel inclusion and perpendicular inclusion of the manifold, i.e.

$$\dot{u} = I_{||}u \qquad\qquad \text{and} \qquad\qquad \dot{r} = I_{\perp}r. \tag{A4.2}$$

With these interpretations, we obtain from equation (3.8) the 2D numerical stress for discrete variables $\boldsymbol{u}, \boldsymbol{r}, \boldsymbol{p}$:

$$\begin{pmatrix} \dot{\boldsymbol{\sigma}} \\ \dot{\boldsymbol{\tau}} \\ \dot{\boldsymbol{v}} \end{pmatrix} = |\varsigma| \begin{pmatrix} -\boldsymbol{\delta}^{-1}\overline{\boldsymbol{\mu}}\Delta^* & -R^{I_{||}\boldsymbol{n}}\tilde{\Xi} & I_{||}\boldsymbol{n}\tilde{\Xi} \\ -R^{I_{||}\boldsymbol{n}}\Xi & & \\ I_{||}\boldsymbol{n}\Xi & & -\boldsymbol{\delta}^{\mu}\Delta^* \end{pmatrix} \begin{pmatrix} I_{||}\boldsymbol{u} \\ I_{\perp}\boldsymbol{r} \\ \boldsymbol{p} \end{pmatrix}. \tag{A4.3}$$

Moreover, we recognize that in 2D, only the two first components of the stresses appear in the finite volume structure, thus by realizing that $I_{||}^*$ (and $I_{\perp}^*$) extracts the two first (and last) components and of a vector we have that:

$$\begin{pmatrix} \boldsymbol{\sigma} \\ \boldsymbol{\tau} \\ \boldsymbol{v} \end{pmatrix} = \begin{pmatrix} I_{||}^*\dot{\boldsymbol{\sigma}} \\ I_{\perp}^*\dot{\boldsymbol{\tau}} \\ \dot{\boldsymbol{v}} \end{pmatrix}. \tag{A4.4}$$

Combining equations (A4.3) and (A4.4), we obtain:

$$\begin{pmatrix} \boldsymbol{\sigma} \\ \boldsymbol{\tau} \\ \boldsymbol{v} \end{pmatrix} = |\boldsymbol{\varsigma}| \begin{pmatrix} -\boldsymbol{\delta}^{-1}\overline{\boldsymbol{\mu}}\Delta^* & -\left(I_{||}^* R^{I_{||}\boldsymbol{n}} I_\perp\right)\tilde{\Xi} & \boldsymbol{n}\tilde{\Xi} \\ -\left(I_\perp^* R^{I_{||}\boldsymbol{n}} I_{||}\right)\Xi & & \\ \boldsymbol{n}\Xi & & -\boldsymbol{\delta}^{\mu}\Delta^* \end{pmatrix} \begin{pmatrix} \boldsymbol{u} \\ \boldsymbol{r} \\ \boldsymbol{p} \end{pmatrix}. \tag{A4.5}$$

A direct calculation shows that:

$$I_{||}^* R^{I_{||}\boldsymbol{n}} I_\perp = I_{||}^* \begin{pmatrix} 0 & 0 & r_2 \\ 0 & 0 & -r_1 \\ -r_2 & r_1 & 0 \end{pmatrix} I_\perp = \begin{pmatrix} r_2 \\ -r_1 \end{pmatrix}, \tag{A4.6}$$

while

$$I_\perp^* R^{I_{||}\boldsymbol{n}} I_{||} = I_{||}^* \begin{pmatrix} 0 & 0 & r_2 \\ 0 & 0 & -r_1 \\ -r_2 & r_1 & 0 \end{pmatrix} I_\perp = \begin{pmatrix} -r_2 & r_1 \end{pmatrix}. \tag{A4.7}$$

This fully specifies the TPSA numerical stresses in 2D.